\documentclass[11pt]{article}%
\usepackage{amsmath}
\usepackage{amsfonts}
\usepackage{amssymb}
\usepackage{graphicx}
\usepackage{a4wide}%
\providecommand{\U}[1]{\protect\rule{.1in}{.1in}}
\providecommand{\U}[1]{\protect\rule{.1in}{.1in}}
\newtheorem{theorem}{Theorem}

\newtheorem{corollary}[theorem]{Corollary}

\newtheorem{example}[theorem]{Example}

\newtheorem{lemma}[theorem]{Lemma}

\newtheorem{proposition}[theorem]{Proposition}
\newtheorem{remark}[theorem]{Remark}

\begin{document}

\title{On constrained optimization problems solved using CDT}
\author{C. Z\u{a}linescu\\Institute of Mathematics \textquotedblleft Octav Mayer\textquotedblright,
Iasi, Romania}
\date{}
\maketitle

\textbf{Abstract} DY Gao together with some of his collaborators applied his
Canonical duality theory (CDT) for solving a class of constrained optimization
problems. Unfortunately, in several papers on this subject there are unclear
statements, not convincing proofs, or even false results. It is our aim in
this work to study rigorously these class of constrained optimization problems
in finite dimensional spaces and to discuss several results published in the
last ten years.

\section{Introduction}

In the preface of the book \emph{Canonical Duality Theory. Advances in
Mechanics and Mathematics, vol 37, Springer, Cham (2017)}, edited by DY Gao, V
Latorre and N Ruan, one says:

\textquotedblleft Canonical duality theory is a breakthrough methodological
theory that can be used not only for modeling complex systems within a unified
framework, but also for solving a large class of challenging problems in
multidisciplinary fields of engineering, mathematics, and sciences. ...

This theory is composed mainly of

(1) a canonical dual transformation, which can be used to formulate perfect
dual problems without duality gap;

(2) a complementary-dual principle, which solved the open problem in finite
elasticity and provides a unified analytical solution form for general
nonconvex/nonsmooth/discrete problems;

(3) a triality theory, which can be used to identify both global and local
optimality conditions and to develop powerful algorithms for solving
challenging problems in complex systems."

\medskip

It is our aim in this work to present rigorously this \textquotedblleft
methodological theory\textquotedblright\ for constrained optimization problems
in finite dimensional spaces. It is not the most general framework, but it
covers all the situations met in the examples provided in DY Gao and his
collaborators' works on constrained optimization problems in finite
dimensions. We also point out some drawbacks and not convincing arguments from
some of those papers.

\section{Preliminaries}

We consider the following minimization problem with equality and inequality constraints

\medskip$(P_{J})$ $\min$ ~$f(x)$ ~~s.t. $x\in X_{J}$, \medskip

\noindent where $J\subset\overline{1,m}$,
\[
X_{J}:=\big\{x\in\mathbb{R}^{n}\mid\left[  \forall j\in J:g_{j}(x)=0\right]
~~\wedge~~\left[  \forall j\in J^{c}:g_{j}(x)\leq0\right]  \big\}
\]
with $J^{c}:=\overline{1,m}\setminus J$, and
\[
f(x):=g_{0}(x):=q_{0}(x)+V_{0}\left(  \Lambda_{0}(x)\right)  ,\quad
g_{j}(x):=q_{j}(x)+V_{j}\left(  \Lambda_{j}(x)\right)  \quad\left(
x\in\mathbb{R}^{n},~~j\in\overline{1,m}\right)  ,
\]
$q_{k}$ and $\Lambda_{k}$ being quadratic functions on $\mathbb{R}^{n}$, and
$V_{k}\in\Gamma_{sc}:=\Gamma_{sc}(\mathbb{R})$ for $k\in\overline{0,m}$; note
that
\begin{equation}
X_{J\cup K}=X_{J}\cap X_{K}\quad\forall J,K\subset\overline{1,m}.
\label{r-xjk}%
\end{equation}

Before giving the precise definition of $\Gamma_{sc}$ we recall some notions
and results from convex analysis we shall use in the sequel.

Having $h:\mathbb{R}^{p}\rightarrow\overline{\mathbb{R}}:=\mathbb{R}%
\cup\{-\infty,+\infty\}$, its domain is $\operatorname*{dom}h:=\{y\in
\mathbb{R}^{p}\mid h(y)<\infty\}$; $h$ is proper if $\operatorname*{dom}%
h\neq\emptyset$ and $h(y)\neq-\infty$ for $y\in\mathbb{R}^{p}$. The Fenchel
conjugate $h^{\ast}:\mathbb{R}^{p}\rightarrow\overline{\mathbb{R}}$ of the
proper function $h$ is defined by
\[
h^{\ast}(\sigma):=\sup\{\left\langle y,\sigma\right\rangle -h(y)\mid
y\in\mathbb{R}^{p}\}=\sup\{\left\langle y,\sigma\right\rangle -h(y)\mid
y\in\operatorname*{dom}h\}\quad(\sigma\in\mathbb{R}^{p}),
\]
while its sub\-differential at $y\in\operatorname*{dom}h$ is
\[
\partial h(y):=\left\{  \sigma\in\mathbb{R}^{p}\mid\left\langle y^{\prime
}-y,\sigma\right\rangle \leq h(y^{\prime})-h(y)~\forall y^{\prime}%
\in\mathbb{R}^{p}\right\}  ,
\]
and $\partial h(y):=\emptyset$ if $y\notin\operatorname*{dom}h$; clearly,
\begin{equation}
h(y)+h^{\ast}(\sigma)\geq\left\langle y,\sigma\right\rangle ~~\wedge~~\left[
\sigma\in\partial h(y)\Longleftrightarrow h(y)+h^{\ast}(\sigma)=\left\langle
y,\sigma\right\rangle \quad\forall(y,\sigma)\in\mathbb{R}^{p}\times
\mathbb{R}^{p}\right]  . \label{r-fen}%
\end{equation}
The class of proper convex lower semi\-continuous (lsc for short) functions
$h:\mathbb{R}^{p}\rightarrow\overline{\mathbb{R}}$ is denoted by
$\Gamma(\mathbb{R}^{p})$. It is well known that for $h\in\Gamma(\mathbb{R}%
^{p})$ one has $h^{\ast}\in\Gamma(\mathbb{R}^{p})$, $(h^{\ast})^{\ast}=h$, and
$\sigma\in\partial h(y)$ iff $y\in\partial h^{\ast}(\sigma)$; moreover,
$\partial h(y)\neq\emptyset$ for every $y\in\operatorname*{ri}%
(\operatorname*{dom}h)$ and $h(\overline{y})=\inf_{y\in\mathbb{R}^{p}}h(y)$
iff $0\in\partial h(\overline{y})$.

We denote by $\Gamma_{sc}(\mathbb{R}^{p})$ the class of those $h\in
\Gamma(\mathbb{R}^{p})$ which are essentially strictly convex and essentially
smooth, that is the class of proper lsc convex functions of Legendre type (see
\cite[Sect.~26]{Roc:72}). For $h\in\Gamma_{sc}(\mathbb{R}^{p})$ we have:
$h^{\ast}\in\Gamma_{sc}(\mathbb{R}^{p})$, $\operatorname*{dom}\partial
h=\operatorname*{int}(\operatorname*{dom}h)$, and $h$ is differentiable on
$\operatorname*{int}(\operatorname*{dom}h)$; moreover, $\nabla
h:\operatorname*{int}(\operatorname*{dom}h)\rightarrow\operatorname*{int}%
(\operatorname*{dom}h^{\ast})$ is bijective and continuous with $\left(
\nabla h\right)  ^{-1}=\nabla h^{\ast}$. Having in view these properties and
(\ref{r-fen}), for $h\in\Gamma_{sc}(\mathbb{R}^{p})$ and $(y,\sigma
)\in\mathbb{R}^{p}\times\mathbb{R}^{p}$ we have that
\begin{align}
h(y)+h^{\ast}(\sigma)=\left\langle y,\sigma\right\rangle  &
\Longleftrightarrow\left[  y\in\operatorname*{int}(\operatorname*{dom}%
h)~\wedge~\sigma=\nabla h(y)\right] \nonumber\\
&  \Longleftrightarrow\left[  \sigma\in\operatorname*{int}(\operatorname*{dom}%
h^{\ast})~\wedge~y=\nabla h^{\ast}(\sigma)\right]  . \label{r-nggs}%
\end{align}
It follows that $\Gamma_{sc}:=\Gamma_{sc}(\mathbb{R})$ is the class of those
$h\in\Gamma(\mathbb{R})$ which are strictly convex and derivable on
$\operatorname*{int}(\operatorname*{dom}h)$, assumed to be nonempty; hence
$h^{\prime}:\operatorname*{int}(\operatorname*{dom}h)\rightarrow
\operatorname*{int}(\operatorname*{dom}h^{\ast})$ is continuous, bijective and
$(h^{\prime})^{-1}=(h^{\ast})^{\prime}$ whenever $h\in\Gamma_{sc}$. The
problem $(P_{\overline{1,m}})$ [resp.\ $(P_{\emptyset})$], denoted by
$(P_{e})$ [resp.\ $(P_{i})$], is a minimization problem with equality
[resp.\ inequality] constraints whose feasible set is $X_{e}:=X_{\overline
{1,m}}$ [resp.\ $X_{i}:=X_{\emptyset}$]. From (\ref{r-xjk}) we get
$X_{e}\subset X_{J}\subset X_{i}$, each inclusion being generally strict for
$J\notin\{\emptyset,\overline{1,m}\}$.

In many examples considered by DY Gao and his collaborators, some functions
$g_{k}$ are quadratic, that is $g_{k}=q_{k}$; we set%
\[
Q:=\{k\in\overline{0,m}\mid g_{k}=q_{k}\},\quad Q_{0}:=Q\setminus
\{0\}=\overline{1,m}\cap Q. \label{r-Q}%
\]
For $k\in Q$ we take $\Lambda_{k}:=0$ and $V_{k}(t):=\tfrac{1}{2}t^{2}$ for
$t\in\mathbb{R}$; then clearly $V_{k}^{\ast}=V_{k}\in\Gamma_{sc}$. To be more
precise, we take%
\[
q_{k}(x):=\tfrac{1}{2}\left\langle x,A_{k}x\right\rangle -\left\langle
b_{k},x\right\rangle +c_{k}~~\wedge~~\Lambda_{k}(x):=\tfrac{1}{2}\left\langle
x,C_{k}x\right\rangle -\left\langle d_{k},x\right\rangle +e_{k}\quad\left(
x\in\mathbb{R}^{n}\right)
\]
with $A_{k},C_{k}\in\mathfrak{S}_{n}$, $b_{k},d_{k}\in\mathbb{R}^{n}$ (seen as
column matrices), and $c_{k},e_{k}\in\mathbb{R}$ for $k\in\overline{0,m}$,
where $\mathfrak{S}_{n}$ denotes the set of $n\times n$ real symmetric
matrices; of course, $c_{0}$ can be taken to be $0$. Clearly, $C_{k}%
=0\in\mathfrak{S}_{n}$, $b_{k}=0\in\mathbb{R}^{n}$ and $c_{k}=0\in\mathbb{R}$
for $k\in Q$. We use also the notations%
\begin{equation}
I_{k}:=\operatorname*{dom}V_{k},\quad I_{k}^{\ast}:=\operatorname*{dom}%
V_{k}^{\ast}\quad(k\in\overline{0,m}),\quad I^{\ast}:={\textstyle\prod
\nolimits_{k=0}^{m}} I_{k}^{\ast}; \label{r-s2}%
\end{equation}
of course, $I_{k}=I_{k}^{\ast}=\mathbb{R}$ for $k\in Q$. In order to simplify
the writing, in the sequel
\[
\lambda_{0}:=\overline{\lambda}_{0}:=1.
\]

To the functions $f$ $(=g_{0})$ and $(g_{j})_{j\in\overline{1,m}}$ we
associate several sets and functions. The Lagrangian $L$ is defined by
\[
L:X\times\mathbb{R}^{m}\rightarrow\mathbb{R},\quad L(x,\lambda):=f(x)+\sum
\nolimits_{j=1}^{m}\lambda_{j}g_{j}(x)=\sum\nolimits_{k=0}^{m}\lambda
_{k}\left[  q_{k}(x)+V_{k}\left(  \Lambda_{k}(x)\right)  \right]
,\label{r-l2}%
\]
where $\lambda:=(\lambda_{1},...,\lambda_{m})^{T}\in\mathbb{R}^{m}$, and
\begin{gather*}
X:=\left\{  x\in\mathbb{R}^{n}\mid\forall k\in\overline{0,m}:\Lambda_{k}%
(x)\in\operatorname*{dom}V_{k}\right\}  ={\textstyle\bigcap\nolimits_{k=0}%
^{m}}\Lambda_{k}^{-1}\left(  \operatorname*{dom}V_{k}\right)  ,\\
X_{0}:=\left\{  x\in\mathbb{R}^{n}\mid\forall k\in\overline{0,m}:\Lambda
_{k}(x)\in\operatorname*{int}(\operatorname*{dom}V_{k})\right\}
\subset\operatorname*{int}X;
\end{gather*}
clearly $X_{0}$ is open and $L$ is differentiable on $X_{0}$. Using Gao's
procedure, we consider the \textquotedblleft extended Lagrangian" $\Xi$
associated to $f$ and $(g_{j})_{j\in\overline{1,m}}$:
\[
\Xi:\mathbb{R}^{n}\times\mathbb{R}^{1+m}\times I^{\ast}\rightarrow
\mathbb{R},\quad\Xi(x,\lambda,\sigma):=\sum\nolimits_{k=0}^{m}\lambda
_{k}\left[  q_{k}(x)+\sigma_{k}\Lambda_{k}(x)-V_{k}^{\ast}(\sigma_{k})\right]
,\label{r-xi2}%
\]
where $I^{\ast}$ is defined in (\ref{r-s2}) and $\sigma:=(\sigma_{0}%
,\sigma_{1},...,\sigma_{m})\in\mathbb{R}\times\mathbb{R}^{m}=\mathbb{R}^{1+m}%
$. Clearly, $\Xi(\cdot,\lambda,\sigma)$ is a quadratic function for every
fixed $(\lambda,\sigma)\in\mathbb{R}^{m}\times I^{\ast}$.

In the sequel we shall use frequently the following sets associated to
$\lambda\in\mathbb{R}^{m}$ and $J\subset\overline{1,m}$:%
\begin{align*}
M_{\neq}(\lambda)  & :=\{j\in\overline{1,m}\mid\lambda_{j}\neq0\},\quad
M_{\neq}^{0}(\lambda):=M_{\neq}(\lambda)\cup\{0\}\\
\Gamma_{J}  & :=\big\{\lambda\in\mathbb{R}^{m}\mid\lambda_{j}\geq0~\forall
j\in J^{c}\big\}\supset\mathbb{R}_{+}^{m},
\end{align*}
respectively; clearly,
\[
\Gamma_{\emptyset}=\mathbb{R}_{+}^{m},\quad\Gamma_{\overline{1,m}}%
=\mathbb{R}^{m},\quad\Gamma_{J\cap K}=\Gamma_{J}\cap\Gamma_{K}\quad\forall
J,K\subset\overline{1,m}.
\]

Taking into account the convexity of the functions $V_{k}$ we obtain useful
relations between $L$ and $\Xi$ in the next result.

\begin{lemma}
\label{lem-xiL}Let $x\in X$ and $J\subset\overline{1,m}$. Then
\begin{equation}
L(x,\lambda)=\sup_{\sigma\in I_{J,Q}}\Xi(x,\lambda,\sigma)\quad\forall
\lambda\in\Gamma_{J\cap Q},\label{r-LXi}%
\end{equation}
where%
\begin{equation}
I_{J,Q}:={\prod\nolimits_{k=0}^{m}}I_{k}^{\ast\ast}\text{~~with~~}I_{k}%
^{\ast\ast}:=\left\{
\begin{array}
[c]{ll}%
\{0\} & \text{if }k\in J\cap Q,\\
I_{k}^{\ast} & \text{if }k\in\overline{0,m}\setminus(J\cap Q),
\end{array}
\right.  \label{r-ikss}%
\end{equation}
and
\[
\sup_{(\lambda,\sigma)\in\Gamma_{J\cap Q}\times I_{J,Q}}\Xi(x,\lambda
,\sigma)=\sup_{\lambda\in\Gamma_{J\cap Q}}L(x,\lambda)=\left\{
\begin{array}
[c]{ll}%
f(x) & \text{if }x\in X_{J\cap Q},\\
\infty & \text{if }x\in X\setminus X_{J\cap Q}.
\end{array}
\right.
\]

\end{lemma}

Proof. Let us set $K:=J\cap Q=J\cap Q_{0}$. It is convenient to observe that
$\Gamma_{K}= {\textstyle\prod\nolimits_{j=1}^{m}} \Gamma_{j}$, where
$\Gamma_{j}:=\mathbb{R}$ for $j\in K$ and $\Gamma_{j}:=\mathbb{R}_{+}$ for
$j\in K^{c}$; moreover, we set $\Gamma_{0}:=\mathbb{R}_{+}$. Take $x\in X$,
$\lambda\in\Gamma_{K}$ and $k\in\overline{0,m}$. Using the fact that
$V_{k}^{\ast\ast}=V_{k}$, we have that
\begin{align*}
g_{k}(x)  &  =q_{k}(x)+V_{k}\left(  \Lambda_{k}(x)\right)  =q_{k}%
(x)+\sup_{\sigma_{k}\in I_{k}^{\ast}}\left[  \sigma_{k}\Lambda_{k}%
(x)-V_{k}^{\ast}(\sigma_{k})\right] \\
&  =\sup_{\sigma_{k}\in I_{k}^{\ast}}\left[  q_{k}(x)+\sigma_{k}\Lambda
_{k}(x)-V_{k}^{\ast}(\sigma_{k})\right]  ,
\end{align*}
whence, because $g_{k}(x)\in\mathbb{R}$,
\begin{equation}
\mu g_{k}(x)=\sup_{\sigma_{k}\in I_{k}^{\ast}}\mu\left[  q_{k}(x)+\sigma
_{k}\Lambda_{k}(x)-V_{k}^{\ast}(\sigma_{k})\right]  \quad\forall\mu
\in\mathbb{R}_{+},~\forall k\in\overline{0,m}. \label{r-lgk}%
\end{equation}
Assume, moreover, that $k\in K$ $(\subset Q_{0}\subset Q)$; then
$g_{k}(x)=q_{k}(x)$, and so
\[
\mu g_{k}(x)=\mu q_{k}(x)=\mu\left[  q_{k}(x)+0\cdot\Lambda_{k}(x)-V_{k}%
^{\ast}(0)\right]  =\sup_{\sigma_{k}\in I_{k}^{\ast\ast}}\mu\left[
q_{k}(x)+\sigma_{k}\Lambda_{k}(x)-V_{k}^{\ast}(\sigma_{k})\right]
\]
for every $\mu\in\mathbb{R}$. Therefore,
\begin{align*}
L(x,\lambda)=  &  \sum_{k\in K}\sup_{\sigma_{k}\in\{0\}}\lambda_{k}\left[
q_{k}(x)+\sigma_{k}\Lambda_{k}(x)-V_{k}^{\ast}(\sigma_{k})\right] \\
&  +\sum_{k\in\overline{0,m}\setminus K}\sup_{\sigma_{k}\in I_{k}^{\ast}%
}\lambda_{k}\left[  q_{k}(x)+\sigma_{k}\Lambda_{k}(x)-V_{k}^{\ast}(\sigma
_{k})\right] \\
=  &  \sum_{k\in\overline{0,m}}\sup_{\sigma_{k}\in I_{k}^{\ast\ast}}%
\lambda_{k}\left[  q_{k}(x)+\sigma_{k}\Lambda_{k}(x)-V_{k}^{\ast}(\sigma
_{k})\right] \\
=  &  \sup_{\sigma\in I_{J,Q}}\sum_{k\in\overline{0,m}}\lambda_{k}\left[
q_{k}(x)+\sigma_{k}\Lambda_{k}(x)-V_{k}^{\ast}(\sigma_{k})\right]
=\sup_{\sigma\in I_{J,Q}}\Xi(x,\lambda,\sigma);
\end{align*}
hence, (\ref{r-LXi}) holds. Using (\ref{r-LXi}) we get
\begin{equation}
\sup_{(\lambda,\sigma)\in\Gamma_{K}\times I_{J,Q}}\Xi(x,\lambda,\sigma
)=\sup_{\lambda\in\Gamma_{K}}\sup_{\sigma\in I_{J,Q}}\Xi(x,\lambda
,\sigma)=\sup_{\lambda\in\Gamma_{K}}L(x,\lambda). \label{r-LXiJ}%
\end{equation}
Since%
\[
\sup_{\lambda\in\mathbb{R}_{+}}\lambda\alpha=\iota_{\mathbb{R}_{-}}%
(\alpha),\quad\sup_{\lambda\in\mathbb{R}}\lambda\alpha=\iota_{\{0\}}(\alpha),
\]
where the indicator function $\iota_{E}:Z\rightarrow\overline{\mathbb{R}}$ of
$E\subset Z$ is defined by $\iota_{E}(z):=0$ for $z\in E$, $\iota
_{E}(z):=+\infty$ for $z\in Z\setminus E$, we get%
\[
\sup_{\lambda\in\Gamma_{K}}L(x,\lambda)=f(x)+\sum_{j\in\overline{1,m}}%
\sup_{\lambda_{j}\in\Gamma_{j}}\lambda_{j}g_{j}(x)=\left\{
\begin{array}
[c]{ll}%
f(x) & \text{if }x\in X_{K},\\
\infty & \text{if }x\in X\setminus X_{K}.
\end{array}
\right.
\]
Using (\ref{r-LXiJ}) and the previous equalities, the conclusion follows.
\hfill$\square$

\medskip

Another useful result in this context is the following.

\begin{lemma}
\label{lem-gj}Let $\overline{x}\in\mathbb{R}^{n}$, $\overline{\sigma}%
\in\mathbb{R}^{m}$ and $k\in\overline{0,m}$. Then
\begin{align*}
\Lambda_{k}(\overline{x})=V_{k}^{\ast\prime}(\overline{\sigma}_{k})  &
\Longleftrightarrow\overline{\sigma}_{k}=V_{k}^{\prime}(\Lambda_{k}%
(\overline{x}))\Longleftrightarrow V_{k}(\Lambda_{k}(\overline{x}%
))+V_{k}^{\ast}(\overline{\sigma}_{k})=\overline{\sigma}_{k}\Lambda
_{k}(\overline{x})\\
&  \Longleftrightarrow g_{k}(\overline{x})=q_{k}(\overline{x})+\overline
{\sigma}_{k}\Lambda_{k}(\overline{x})-V_{k}^{\ast}(\overline{\sigma}%
_{k})\label{r-gkqk}\\
&  \Longrightarrow\left[  \overline{\sigma}_{k}\in\operatorname*{int}%
(\operatorname*{dom} V_{k}^{\ast})~~\wedge~~\Lambda_{k}(\overline{x}%
)\in\operatorname*{int} (\operatorname*{dom}V_{k})\right]  .
\end{align*}

In particular, for $k\in Q$, $\Lambda_{k}(\overline{x})=V_{k}^{\ast\prime
}(\overline{\sigma}_{k})$ if and only if $\overline{\sigma}_{k}=0$.
\end{lemma}

Proof. Because $V_{k}\in\Gamma_{sc}$, (\ref{r-nggs}) holds. Since
$g_{k}(\overline{x})=q_{k}(\overline{x})+V_{k}(\Lambda_{k}(\overline{x}))$, we
obtain that $g_{k}(\overline{x})=q_{k}(\overline{x})+\overline{\sigma}%
_{k}\Lambda_{k}(\overline{x})-V_{k}^{\ast}(\overline{\sigma}_{k})$ if and only
$V_{k}(\Lambda_{k}(\overline{x}))=\overline{\sigma}_{k}\Lambda_{k}%
(\overline{x})-V_{k}^{\ast}(\overline{\sigma}_{k})$, and so the conclusion
follows. The case $k\in Q$ follows immediately. \hfill$\square$

\medskip

\begin{corollary}
\label{cor-LXi}Let $(\overline{x},\overline{\lambda},\overline{\sigma})\in
X\times\mathbb{R}^{m}\times I^{\ast}$ with $\overline{\sigma}_{k}=0$ for $k\in
Q$. If
\begin{equation}
\forall k\in M_{\neq}^{0}(\overline{\lambda})\setminus Q:\left[  \Lambda
_{k}(\overline{x})\in\operatorname*{int}(\operatorname*{dom}V_{k}%
)~\wedge~\overline{\sigma}_{k}\in\operatorname*{int}(\operatorname*{dom}%
V_{k}^{\ast})~\wedge~\overline{\sigma}_{k}=V_{k}^{\prime}\left(  \Lambda
_{k}(\overline{x})\right)  \right]  , \label{r-ml0}%
\end{equation}
then $L(\overline{x},\overline{\lambda})=\Xi(\overline{x},\overline{\lambda
},\overline{\sigma})$. Conversely, if $L(\overline{x},\overline{\lambda}%
)=\Xi(\overline{x},\overline{\lambda},\overline{\sigma})$ and $\overline
{\lambda}\in\Gamma_{Q_{0}}$, then (\ref{r-ml0}) holds.
\end{corollary}

Proof. Assume first that (\ref{r-ml0}) holds. Using Lemma \ref{lem-gj} we
obtain that $V_{k}\left(  \Lambda_{k}(\overline{x})\right)  =\overline{\sigma
}_{k}\Lambda_{k}(\overline{x})-V_{k}^{\ast}(\overline{\sigma}_{k})$ for $k\in
M_{\neq}^{0}(\overline{\lambda})$, and so
\[
\overline{\lambda}_{k}\left[  q_{k}(\overline{x})+V_{k}\left(  \Lambda
_{k}(\overline{x})\right)  \right]  =\overline{\lambda}_{k}\left[
q_{k}(\overline{x})+\overline{\sigma}_{k}\Lambda_{k}(\overline{x})-V_{k}%
^{\ast}(\overline{\sigma}_{k})\right]  \quad\forall k\in\overline
{0,m}\label{r-kLXi}%
\]
because $\overline{\lambda}_{k}=0$ for $k\not \in \overline{0,m}\setminus
M_{\neq}^{0}(\overline{\lambda})$. Then the equality $L(\overline{x}%
,\overline{\lambda})=\Xi(\overline{x},\overline{\lambda},\overline{\sigma})$
follows from de definitions of $L$ and $\Xi$.

Conversely, assume that $L(\overline{x},\overline{\lambda})=\Xi(\overline
{x},\overline{\lambda},\overline{\sigma})$ and $\overline{\lambda}\in
\Gamma_{Q_{0}}$; hence $\overline{\lambda}_{k}\geq0$ for all $k\in Q_{0}^{c}$.
Clearly, $g_{k}(\overline{x})=q_{k}(\overline{x})=q_{k}(\overline
{x})+\overline{\sigma}_{k}\Lambda_{k}(\overline{x})-V_{k}^{\ast}%
(\overline{\sigma}_{k})$ for $k\in Q$. Because $g_{k}(\overline{x}%
)=q_{k}(\overline{x})+V_{k}\left(  \Lambda_{k}(\overline{x})\right)  \geq
q_{k}(\overline{x})+\overline{\sigma}_{k}\Lambda_{k}(\overline{x})-V_{k}%
^{\ast}(\overline{\sigma}_{k})$ and $\overline{\lambda}_{k}\geq0$ for all
$k\in\{0\}\cup Q_{0}^{c}\supset\overline{0,m}\setminus Q$, from $L(\overline
{x},\overline{\lambda})=\Xi(\overline{x},\overline{\lambda},\overline{\sigma
})$ we obtain that
\[
\overline{\lambda}_{k}\left[  q_{k}(\overline{x})+V_{k}\left(  \Lambda
_{k}(\overline{x})\right)  \right]  =\overline{\lambda}_{k}\left[
q_{k}(\overline{x})+\overline{\sigma}_{k}\Lambda_{k}(\overline{x})-V_{k}%
^{\ast}(\overline{\sigma}_{k})\right]  \quad\forall k\in\overline
{0,m}\setminus Q;
\]
hence $g_{k}(\overline{x})=q_{k}(\overline{x})+\overline{\sigma}_{k}%
\Lambda_{k}(\overline{x})-V_{k}^{\ast}(\overline{\sigma}_{k})$, that is
$V_{k}\left(  \Lambda_{k}(\overline{x})\right)  =\overline{\sigma}_{k}%
\Lambda_{k}(\overline{x})-V_{k}^{\ast}(\overline{\sigma}_{k})$, for $k\in
M_{\neq}^{0}(\overline{\lambda})\setminus Q$. Using (\ref{r-nggs}) we obtain
that (\ref{r-ml0}) is verified. \hfill$\square$

\medskip

Let us consider $G:\mathbb{R}^{m}\times\mathbb{R}^{1+m}\rightarrow
\mathfrak{S}_{n}$, $F:\mathbb{R}^{m}\times\mathbb{R}^{1+m}\rightarrow
\mathbb{R}^{n}$, $E:\mathbb{R}^{m}\times\mathbb{R}^{1+m}\rightarrow\mathbb{R}$
defined by
\[
G(\lambda,\sigma):=\sum_{k=0}^{m}\lambda_{k}(A_{k}+\sigma_{k}C_{k}%
),~~F(\lambda,\sigma):=\sum_{k=0}^{m}\lambda_{k}(b_{k}+\sigma_{k}%
d_{k}),~~E(\lambda,\sigma):=\sum_{k=0}^{m}\lambda_{k}(c_{k}+\sigma_{k}e_{k}).
\label{r-cls}%
\]
Hence, for $(\lambda,\sigma)\in\mathbb{R}^{m}\times I^{\ast}$ we have that
\[
\Xi(x,\lambda,\sigma)=\tfrac{1}{2}\left\langle x,G(\lambda,\sigma
)x\right\rangle -\left\langle F(\lambda,\sigma),x\right\rangle +E(\lambda
,\sigma)-\sum\nolimits_{k=0}^{m}\lambda_{k}V_{k}^{\ast}(\sigma_{k}).
\label{r-x2}%
\]

\begin{remark}
\label{rem-afin}Note that $G$, $F$ and $E$ do not depend on $\sigma_{k}$ for
$k\in Q$. Moreover, $G$, $F$ and $E$ are affine functions when $\overline
{1,m}\subset Q$, that is $Q_{0}=\overline{1,m}$.
\end{remark}

For $(\lambda,\sigma)\in\mathbb{R}^{m}\times I^{\ast}$ we have that
\begin{gather}
\nabla_{x}\Xi(x,\lambda,\sigma)=G(\lambda,\sigma)x-F(\lambda,\sigma
),\quad\nabla_{xx}^{2}\Xi(x,\lambda,\sigma)=G(\lambda,\sigma), \label{r-d1xx2}%
\\
\nabla_{\lambda}\Xi(x,\lambda,\sigma)=\left(  q_{1}(x)+\sigma_{1}\Lambda
_{1}(x)-V_{1}^{\ast}(\sigma_{1}),...,q_{m}(x)+\sigma_{m}\Lambda_{m}%
(x)-V_{m}^{\ast}(\sigma_{m})\right)  ^{T}, \label{r-d1mx2}%
\end{gather}
while for $(x,\lambda,\sigma)\in\mathbb{R}^{n}\times\mathbb{R}^{m}%
\times\operatorname*{int}I^{\ast}$ we have that
\begin{equation}
\nabla_{\sigma}\Xi(x,\lambda,\sigma)=\left(  \lambda_{0}\left[  \Lambda
_{0}(x)-V_{0}^{\ast\prime}(\sigma_{0})\right]  ,\lambda_{1}\left[  \Lambda
_{1}(x)-V_{1}^{\ast\prime}(\sigma_{1})\right]  ,...,\lambda_{m}\left[
\Lambda_{m}(x)-V_{m}^{\ast\prime}(\sigma_{m})\right]  \right)  ^{T}.
\label{r-d1sx2}%
\end{equation}

Other relations between $L$ and $\Xi$ are provided in the next result.

\begin{lemma}
\label{lem-nXiL}Let $(\overline{x},\overline{\lambda},\overline{\sigma})\in
X_{0}\times\mathbb{R}^{m}\times\operatorname*{int}I^{\ast}$ be such that
$\nabla_{\sigma}\Xi(\overline{x},\overline{\lambda},\overline{\sigma})=0$ and
$\overline{\sigma}_{k}=0$ for $k\in Q$. Then $L(\overline{x},\overline
{\lambda})=\Xi(\overline{x},\overline{\lambda},\overline{\sigma})$ and
$\nabla_{x}L(\overline{x},\overline{\lambda})=\nabla_{x}\Xi(\overline
{x},\overline{\lambda},\overline{\sigma})$. Moreover, for $j\in\overline{1,m}%
$, $\frac{\partial L}{\partial\lambda_{j}}(\overline{x},\overline{\lambda
})\geq\frac{\partial\Xi}{\partial\lambda_{j}}(\overline{x},\overline{\lambda
},\overline{\sigma})$, with equality if $j\in M_{\neq}(\overline{\lambda})\cup
Q_{0}$; in particular $\nabla_{\lambda}L(\overline{x},\overline{\lambda
})=\nabla_{\lambda}\Xi(\overline{x},\overline{\lambda},\overline{\sigma})$ if
$M_{\neq}(\overline{\lambda})\supset Q_{0}^{c}$ $(=\overline{1,m}\setminus Q)$.
\end{lemma}

Proof. For $k\in M_{\neq}^{0}(\overline{\lambda})$ we have that $\overline
{\lambda}_{k}\neq0$; using (\ref{r-d1sx2}) and Lemma \ref{lem-gj} we get
$\Lambda_{k}(\overline{x})-V_{k}^{\ast\prime}(\overline{\sigma}_{k})=0$, and
so $\overline{\sigma}_{k}\in\operatorname*{int}(\operatorname*{dom}V_{k}%
^{\ast})$, $\Lambda_{k}(\overline{x})\in\operatorname*{int}%
(\operatorname*{dom}V_{k})$, $\overline{\sigma}_{k}=V_{k}^{\prime}(\Lambda
_{k}(\overline{\overline{x}}))$ for $k\in M_{\neq}^{0}(\overline{\lambda})$.
Hence (\ref{r-ml0}) is verified, and so $L(\overline{x},\overline{\lambda
})=\Xi(\overline{x},\overline{\lambda},\overline{\sigma})$ by Corollary
\ref{cor-LXi}. Moreover,%
\begin{align*}
\nabla_{x}L(\overline{x},\overline{\lambda})  &  =\sum\nolimits_{k\in
\overline{0,m}}\overline{\lambda}_{k}\left[  A_{k}\overline{x}-b_{k}%
+V_{k}^{\prime}(\Lambda_{k}(\overline{x}))(C_{k}\overline{x}-d_{k})\right] \\
&  =\sum\nolimits_{k\in M_{\neq}^{0}(\overline{\lambda})}\overline{\lambda
}_{k}\left[  A_{k}\overline{x}-b_{k}+V_{k}^{\prime}(\Lambda_{k}(\overline
{x}))(C_{k}\overline{x}-d_{k})\right]  ,\\
\nabla_{x}\Xi(\overline{x},\overline{\lambda},\overline{\sigma})  &
=\sum\nolimits_{k\in\overline{0,m}}\overline{\lambda}_{k}\left[
A_{k}\overline{x}-b_{k}+\overline{\sigma}_{k}(C_{k}\overline{x}-d_{k})\right]
\\
&  =\sum\nolimits_{k\in M_{\neq}^{0}(\overline{\lambda})}\overline{\lambda
}_{k}\left[  A_{k}\overline{x}-b_{k}+\overline{\sigma}_{k}(C_{k}\overline
{x}-d_{k})\right]  ,
\end{align*}
and so $\nabla_{x}L(\overline{x},\overline{\lambda})=\nabla_{x}\Xi
(\overline{x},\overline{\lambda},\overline{\sigma})$. Clearly, from the
definitions of $L$, $\Xi$ and the inequality in (\ref{r-fen}), we have that
\[
\frac{\partial L}{\partial\lambda_{j}}(\overline{x},\overline{\lambda}%
)=g_{j}(\overline{x})=q_{j}(\overline{x})+\overline{\sigma}_{j}V_{j}%
(\Lambda_{j}(\overline{x}))\geq q_{j}(\overline{x})+\overline{\sigma}%
_{j}\Lambda_{j}(\overline{x})-V_{j}^{\ast}(\overline{\sigma}_{j}%
)=\frac{\partial\Xi}{\partial\lambda_{j}}(\overline{x},\overline{\lambda
},\overline{\sigma}).
\]
Using Lemma \ref{lem-gj} we obtain that $g_{k}(\overline{x})=q_{k}%
(\overline{x})+\overline{\sigma}_{k}\Lambda_{k}(\overline{x})-V_{k}^{\ast
}(\overline{\sigma}_{k})$ for $k\in M_{\neq}(\overline{\lambda})\cup Q$ and so
$\frac{\partial L}{\partial\lambda_{j}}(\overline{x},\overline{\lambda}%
)=\frac{\partial\Xi}{\partial\lambda_{j}}(\overline{x},\overline{\lambda
},\overline{\sigma})$ for $j\in M_{\neq}(\overline{\lambda})\cup Q_{0}$.
\hfill$\square$

\medskip

We consider also the sets
\begin{gather*}
T_{Q}:=\left\{  (\lambda,\sigma)\in\mathbb{R}^{m}\times I^{\ast}\mid\det
G(\lambda,\sigma)\neq\emptyset\wedge\lbrack\forall k\in Q:\sigma
_{k}=0]\right\}  ,\\
T_{Q,\operatorname{col}}:=\left\{  (\lambda,\sigma)\in\mathbb{R}^{m}\times
I^{\ast}\mid F(\lambda,\sigma)\in\operatorname{Im}G(\lambda,\sigma
)\wedge\lbrack\forall k\in Q:\sigma_{k}=0]\right\}  \supseteq T_{Q},\\
T_{Q}^{J+}:=\left\{  (\lambda,\sigma)\in T_{Q}\mid\lambda\in\Gamma_{J\cap
Q},~G(\lambda,\sigma)\succ0\right\}  ,\\
T_{Q,\operatorname{col}}^{J+}:=\left\{  (\lambda,\sigma)\in
T_{Q,\operatorname{col}}\mid\lambda\in\Gamma_{J\cap Q},~G(\lambda
,\sigma)\succeq0\right\}  \supseteq T_{Q}^{J+},
\end{gather*}
as well as the sets%
\[
T:=T_{\emptyset},\quad T_{\operatorname{col}}:=T_{\emptyset,\operatorname{col}%
},\quad T^{+}:=T_{\emptyset}^{\emptyset+},\quad T_{\operatorname{col}}%
^{+}:=T_{\emptyset,\operatorname{col}}^{\emptyset+};
\]
in general $T_{Q}^{J+}$ and $T_{Q,\operatorname{col}}^{J+}$ are not convex,
unlike their corresponding sets $Y^{+}$, $Y_{\operatorname{col}}^{+}$ and
$S^{+}$, $S_{\operatorname{col}}^{+}$ from \cite{Zal:18b} and \cite{Zal:18c},
respectively. However, taking into account Remark \ref{rem-afin},
$T_{\operatorname{col}}$, $T_{Q}^{J+}$ and $T_{Q,\operatorname{col}}^{J+}$ are
convex whenever $Q_{0}=\overline{1,m}$. In the present context it is natural
(in fact necessary) to take $\lambda\in\Gamma_{Q_{0}}$. As in \cite{Zal:18b}
and \cite{Zal:18c}, we consider the (dual objective) function
\[
D:T_{\operatorname{col}}\rightarrow\mathbb{R},\quad D(\lambda,\sigma
):=\Xi(x,\lambda,\sigma)\text{ with }G(\lambda,\sigma)x=F(\lambda,\sigma);
\]
$D$ is well defined by \cite[Lem.~1 (ii)]{Zal:18b}. Consider
\begin{equation}
\xi:T\rightarrow\mathbb{R}^{n},\quad\xi(\lambda,\sigma):=G(\lambda
,\sigma)^{-1}F(\lambda,\sigma). \label{r-xls}%
\end{equation}
For $(\lambda,\sigma)\in T$ we obtain that%
\begin{align}
D(\lambda,\sigma)  &  =\Xi(G(\lambda,\sigma)^{-1}F(\lambda,\sigma
),\lambda,\sigma)=\Xi(\xi(\lambda,\sigma),\lambda,\sigma)\nonumber\\
&  =-\tfrac{1}{2}\left\langle F(\lambda,\sigma),G(\lambda,\sigma
)^{-1}F(\lambda,\sigma)\right\rangle +E(\lambda,\sigma)-\sum\nolimits_{k=0}%
^{m}\lambda_{k}V_{k}^{\ast}(\sigma_{k}). \label{r-p2d}%
\end{align}
Taking into account the second formula in (\ref{r-d1xx2}), we have that
$\Xi(\cdot,\lambda,\sigma)$ is [strictly] convex for $(\lambda,\sigma)\in
T_{\operatorname{col}}^{+}$ $[(\lambda,\sigma)\in T^{+}]$, and so
\begin{equation}
D(\lambda,\sigma)=\min_{x\in\mathbb{R}^{n}}\Xi(x,\lambda,\sigma)\quad
\forall(\lambda,\sigma)\in T_{\operatorname{col}}\text{ such that }%
G(\lambda,\sigma)\succeq0, \label{r-pd2}%
\end{equation}
the minimum being attained uniquely at $\xi(\lambda,\sigma)$ when, moreover,
$G(\lambda,\sigma)\succ0$.

\begin{proposition}
\label{p-perfdual}Let $(\overline{x},\overline{\lambda},\overline{\sigma}%
)\in\mathbb{R}^{n}\times\mathbb{R}^{m}\times I^{\ast}$ be such that
$\nabla_{x}\Xi(\overline{x},\overline{\lambda},\overline{\sigma})=0$,
$\frac{\partial\Xi}{\partial\sigma_{0}}(\overline{x},\overline{\lambda
},\overline{\sigma})=0$, and $\left\langle \overline{\lambda},\nabla_{\lambda
}\Xi(\overline{x},\overline{\lambda},\overline{\sigma})\right\rangle =0$. Then
$(\overline{\lambda},\overline{\sigma})\in T_{\operatorname{col}}$ and%
\begin{equation}
f(\overline{x})=\Xi(\overline{x},\overline{\lambda},\overline{\sigma
})=D(\overline{\lambda},\overline{\sigma}). \label{r-px2pd}%
\end{equation}

\end{proposition}

Proof. Because $\nabla_{x}\Xi(\overline{x},\overline{\lambda},\overline
{\sigma})=0$, $(\overline{\lambda},\overline{\sigma})\in T_{\operatorname{col}%
}$ and the second equality in (\ref{r-px2pd}) holds by the definition of $D$.
Since $\Lambda_{0}(\overline{x})-V_{0}^{\ast\prime}(\overline{\sigma}%
_{0})=\frac{\partial\Xi}{\partial\sigma_{0}}(\overline{x},\overline{\lambda
},\overline{\sigma})=0$, we have that $V_{0}\left(  \Lambda_{0}(\overline
{x})\right)  =\overline{\sigma}_{0}\Lambda_{0}(\overline{x})-V_{0}^{\ast
}(\overline{\sigma}_{0})$ by Lemma \ref{lem-gj} for $k:=0$. Therefore,
\[
f(\overline{x})=q_{0}(\overline{x})+V_{0}\left(  \Lambda_{0}(\overline
{x})\right)  =q_{0}(\overline{x})+\overline{\sigma}_{0}\Lambda_{0}%
(\overline{x})-V^{\ast}(\overline{\sigma}_{0}),
\]
whence
\[
\Xi(\overline{x},\overline{\lambda},\overline{\sigma})=\overline{\lambda}%
_{0}\left[  q_{0}(x)+\sigma_{0}\Lambda_{0}(x)-V_{0}^{\ast}(\sigma_{0})\right]
+\left\langle \overline{\lambda},\nabla_{\lambda}\Xi(\overline{x}%
,\overline{\lambda},\overline{\sigma})\right\rangle =f(\overline{x}).
\]
Hence the first equality in (\ref{r-px2pd}) holds, too. \hfill$\square$

\medskip

Formula (\ref{r-px2pd}) is related to the so-called \textquotedblleft
complimentary-dual principle\textquotedblright\ (see \cite[p.~NP11]%
{GaoRuaLat:16}, \cite[p.~13]{GaoRuaLat:17}) and sometimes is called the
\textquotedblleft perfect duality formula\textquotedblright.

Observe that $T\cap(\mathbb{R}^{m}\times\operatorname*{int}I^{\ast}%
)\subset\operatorname*{int}T$, and for any $\overline{\sigma}\in I^{\ast}$ we
have that the set $\{\lambda\in\mathbb{R}^{m}\mid(\lambda,\overline{\sigma
})\in T\}$ is open. Similarly to the computation of $\frac{\partial
D(\lambda)}{\partial\lambda_{j}}$ in \cite[p.~5]{Zal:18b}, using the
expression of $D(\lambda,\sigma)$ in (\ref{r-p2d}), we get
\begin{align}
\frac{\partial D(\lambda,\sigma)}{\partial\lambda_{j}}  &  =\tfrac{1}%
{2}\left\langle \xi(\lambda,\sigma),(A_{j}+\sigma_{j}C_{j})\xi(\lambda
,\sigma)\right\rangle -\left\langle b_{j}+\sigma_{j}d_{j},\xi(\lambda
,\sigma)\right\rangle +c_{j}+\sigma_{j}e_{j}-V_{j}^{\ast}(\sigma
_{j})\nonumber\\
&  =q_{j}\left(  \xi(\lambda,\sigma)\right)  +\sigma_{j}\Lambda_{j}\left(
\xi(\lambda,\sigma)\right)  -V_{j}^{\ast}(\sigma_{j})\quad\forall
j\in\overline{1,m},~~\forall(\lambda,\sigma)\in T, \label{r-d1p2m}%
\end{align}
and
\begin{align}
\frac{\partial D(\lambda,\sigma)}{\partial\sigma_{k}}  &  =\lambda_{k}\left[
\tfrac{1}{2}\left\langle \xi(\lambda,\sigma),C_{k}\xi(\lambda,\sigma
)\right\rangle -\left\langle d_{k},\xi(\lambda,\sigma)\right\rangle
+e_{k}-V_{k}^{\ast\prime}(\sigma_{k})\right] \nonumber\\
&  =\lambda_{k}\left[  \Lambda_{k}\left(  \xi(\lambda,\sigma)\right)
-V_{k}^{\ast\prime}(\sigma_{k})\right]  \quad\forall k\in\overline
{0,m},~~\forall(\lambda,\sigma)\in T\cap(\mathbb{R}^{m}\times
\operatorname*{int}I^{\ast}). \label{r-d1p2s}%
\end{align}

\begin{lemma}
\label{lem1}Let $(\overline{\lambda},\overline{\sigma})\in\left(
\mathbb{R}^{m}\times\operatorname*{int}I^{\ast}\right)  \cap T$ and set
$\overline{x}:=\xi(\overline{\lambda},\overline{\sigma})$. Then
\[
\nabla_{x}\Xi(\overline{x},\overline{\lambda},\overline{\sigma})=0~~\wedge
~~\nabla_{\lambda}\Xi(\overline{x},\overline{\lambda},\overline{\sigma
})=\nabla_{\lambda}D(\overline{\lambda},\overline{\sigma})~~\wedge
~~\nabla_{\sigma}\Xi(\overline{x},\overline{\lambda},\overline{\sigma}%
)=\nabla_{\sigma}D(\overline{\lambda},\overline{\sigma}).
\]
In particular $(\overline{x},\overline{\lambda},\overline{\sigma})$ is a
critical point of $\Xi$ if and only if $(\overline{\lambda},\overline{\sigma
})$ is a critical point of $D$.
\end{lemma}

Proof. Using (\ref{r-d1xx2}) we get $\nabla_{x}\Xi(\overline{x},\overline
{\lambda},\overline{\sigma})=0$. From (\ref{r-d1p2m}) and (\ref{r-d1mx2}) for
$j\in\overline{1,m}$ we get
\[
\frac{\partial D}{\partial\lambda_{j}}(\overline{\lambda},\overline{\sigma
})=q_{j}(\overline{x})+\overline{\sigma}_{j}\Lambda_{j}(\overline{x}%
)-V_{j}^{\ast}(\overline{\sigma}_{j})=\frac{\partial\Xi}{\partial\lambda_{j}%
}(\overline{x},\overline{\lambda},\overline{\sigma}),
\]
while from (\ref{r-d1p2s}) and (\ref{r-d1sx2}) for $k\in\overline{0,m}$ we get%
\[
\frac{\partial D}{\partial\sigma_{k}}(\overline{\lambda},\overline{\sigma
})=\overline{\lambda}_{k}\left[  \Lambda_{k}(\overline{x})-V_{k}^{\ast\prime
}(\overline{\sigma}_{k})\right]  =\frac{\partial\Xi}{\partial\sigma_{k}%
}(\overline{x},\overline{\lambda},\overline{\sigma}).
\]
The conclusion follows. \hfill$\square$

\medskip

Similarly to \cite{Zal:18b}, we say that $(\overline{x},\overline{\lambda})\in
X_{0}\times\mathbb{R}^{m}$ is a $J$-LKKT point of $L$ if $\nabla
_{x}L(\overline{x},\overline{\lambda})=0$ and
\[
\big[\forall j\in J^{c}:\overline{\lambda}_{j}\geq0~\wedge~\frac{\partial
L}{\partial\lambda_{j}}(\overline{x},\overline{\lambda})\leq0~\wedge
~\overline{\lambda}_{j}\frac{\partial L}{\partial\lambda_{j}}(\overline
{x},\overline{\lambda})=0\big]~\wedge~\big[\forall j\in J:\frac{\partial
L}{\partial\lambda_{j}}(\overline{x},\overline{\lambda})=0\big],
\]
or, equivalently,
\[
\overline{x}\in X_{J}~~\wedge~~\overline{\lambda}\in\Gamma_{J}~~\wedge
~~\left[  \forall j\in J^{c}:\overline{\lambda}_{j}g_{j}(\overline
{x})=0\right]  ;
\]
moreover, we say that $\overline{x}\in X_{0}$ is a $J$-LKKT point of $(P_{J})$
if there exists $\overline{\lambda}\in\mathbb{R}^{m}$ such that $(\overline
{x},\overline{\lambda})$ is a $J$-LKKT point of $L$. Inspired by these
notions, we say that $(\overline{x},\overline{\lambda},\overline{\sigma}%
)\in\mathbb{R}^{n}\times\mathbb{R}^{m}\times\operatorname*{int}I^{\ast}$ is a
$J$-LKKT point of $\Xi$ if $\nabla_{x}\Xi(\overline{x},\overline{\lambda
},\overline{\sigma})=0$, $\nabla_{\sigma}\Xi(\overline{x},\overline{\lambda
},\overline{\sigma})=0$ and
\[
\big[\forall j\in J^{c}:\overline{\lambda}_{j}\geq0~\wedge~\frac{\partial\Xi
}{\partial\lambda_{j}}(\overline{x},\overline{\lambda},\overline{\sigma}%
)\leq0~\wedge~\overline{\lambda}_{j}\frac{\partial\Xi}{\partial\lambda_{j}%
}(\overline{x},\overline{\lambda},\overline{\sigma})=0\big]~\wedge
~\big[\forall j\in J:\frac{\partial\Xi}{\partial\lambda_{j}}(\overline
{x},\overline{\lambda},\overline{\sigma})=0\big], \label{r-kkt-x2}%
\]
and $(\overline{\lambda},\overline{\sigma})\in\left(  \mathbb{R}^{m}%
\times\operatorname*{int}I^{\ast}\right)  \cap T$ is a $J$-LKKT point of $D$
if $\nabla_{\sigma}D(\overline{\lambda},\overline{\sigma})=0$ and
\[
\big[\forall j\in J^{c}:\overline{\lambda}_{j}\geq0~\wedge~\frac{\partial
D}{\partial\lambda_{j}}(\overline{\lambda},\overline{\sigma})\leq
0~\wedge~\overline{\lambda}_{j}\frac{\partial D}{\partial\lambda_{j}%
}(\overline{\lambda},\overline{\sigma})=0\big]~\wedge~\big[\forall j\in
J:\frac{\partial D}{\partial\lambda_{j}}(\overline{\lambda},\overline{\sigma
})=0\big].
\]

In the case in which $J=\emptyset$ we obtain the notions of KKT points for
$\Xi$ and $D$. So, $(\overline{x},\overline{\lambda},\overline{\sigma}%
)\in\mathbb{R}^{n}\times\mathbb{R}^{m}\times\operatorname*{int}I^{\ast}$ is a
KKT point of $\Xi$ if $\nabla_{x}\Xi(\overline{x},\overline{\lambda}%
,\overline{\sigma})=0$, $\nabla_{\sigma}\Xi(\overline{x},\overline{\lambda
},\overline{\sigma})=0$ and%
\begin{equation}
\overline{\lambda}\in\mathbb{R}_{+}^{m}~~\wedge~~\nabla_{\lambda}\Xi
(\overline{x},\overline{\lambda},\overline{\sigma})\in\mathbb{R}_{-}%
^{m}~~\wedge~~\,\left\langle \overline{\lambda},\nabla_{\lambda}\Xi
(\overline{x},\overline{\lambda},\overline{\sigma})\right\rangle =0,
\label{r-k2t-x2i}%
\end{equation}
and $(\overline{\lambda},\overline{\sigma})\in\mathbb{R}^{m}\times
\operatorname*{int}I^{\ast}$ is a KKT point of $D$ if $\nabla_{\sigma
}D(\overline{\lambda},\overline{\sigma})=0$ and%
\[
\overline{\lambda}\in\mathbb{R}_{+}^{m}~~\wedge~~\nabla_{\lambda}%
D(\overline{\lambda},\overline{\sigma})\in\mathbb{R}_{-}^{m}~~\wedge
~~\,\big\langle\overline{\lambda},\nabla_{\lambda}D(\overline{\lambda
},\overline{\sigma})\big\rangle=0. \label{r-k2t-pdi}%
\]

\begin{remark}
\label{rem-kkt}The definition of a KKT point for $\Xi$ is suggested in the
proof of \cite[Th.\ 3]{RuaGao:17} (the same as that of \cite[Th.\ 3]%
{RuaGao:16}). Observe that $(\overline{x},\overline{\lambda},\overline{\sigma
})$ verifying the conditions in (\ref{r-k2t-x2i}) is called critical point of
$\Xi$ in \cite[p.~477]{GaoRuaShe:09}.
\end{remark}

\begin{corollary}
\label{c-lkkt}Let $(\overline{\lambda},\overline{\sigma})\in\left(
\mathbb{R}^{m}\times\operatorname*{int}I^{\ast}\right)  \cap T$.

\emph{(i)} If $\overline{x}:=\xi(\overline{\lambda},\overline{\sigma})$, then
$(\overline{x},\overline{\lambda},\overline{\sigma})$ is a $J$-LKKT point of
$\Xi$ if and only if $(\overline{\lambda},\overline{\sigma})\ $is a $J$-LKKT
point of $D$.

\emph{(ii)} If $M_{\neq}(\overline{\lambda})=\overline{1,m}$, then
$(\overline{x},\overline{\lambda},\overline{\sigma})$ is a $J$-LKKT point of
$\Xi$ if and only if $(\overline{x},\overline{\lambda},\overline{\sigma})$ is
a critical point of $\Xi$, if and only if $\overline{x}=\xi(\overline{\lambda
},\overline{\sigma})$ and $(\overline{\lambda},\overline{\sigma})\ $is a
critical point of $D$.
\end{corollary}

Proof. (i) is immediate from Lemma \ref{lem1}, while (ii) is an obvious
consequence of (i) and the definitions of the corresponding notions.
\hfill$\square$

\begin{remark}
\label{rem-skQ}Taking into account Remark \ref{rem-afin}, as well as
(\ref{r-d1xx2}), (\ref{r-xls}) and Lemma \ref{lem1}, the functions $\nabla
_{x}\Xi$, $\xi$, $\nabla_{\sigma}D$ do not depend on $\sigma_{k}$ for $k\in
Q$. Consequently, if $(\overline{x},\overline{\lambda},\overline{\sigma})$ is
a $J$-LKKT point of $\Xi$ then $\overline{\sigma}_{k}=0$ for $k\in Q\cap
M_{\neq}(\overline{\lambda})$, and $(\overline{x},\overline{\lambda}%
,\tilde{\sigma})$ is also a $J$-LKKT point of $\Xi$, where $\tilde{\sigma}%
_{k}:=0$ for $k\in Q$ and $\tilde{\sigma}_{k}:=\overline{\sigma}_{k}$ for
$k\in\overline{0,m}\setminus Q$. Conversely, taking into account that
$\nabla_{\sigma}D$ does not depend on $\sigma_{k}$ for $k\in Q$, if
$(\overline{\lambda},\overline{\sigma})\in T$ is a $J$-LKKT point of $D$ then
$(\overline{\lambda},\tilde{\sigma})$ is also a $J$-LKKT point of $D$, where
$\tilde{\sigma}_{k}:=0$ for $k\in Q$ and $\tilde{\sigma}_{k}:=\overline
{\sigma}_{k}$ for $k\in\overline{0,m}\setminus Q$.
\end{remark}

Having in view the previous remark, without loss of generality, in the sequel
we shall assume that $\overline{\sigma}_{k}=0$ for $k\in Q$ when
$(\overline{x},\overline{\lambda},\overline{\sigma})\in\mathbb{R}^{n}%
\times\mathbb{R}^{m}\times\operatorname*{int}I^{\ast}$ is a $J$-LKKT point of
$\Xi$, or $(\overline{\lambda},\overline{\sigma})\in T$ is a $J$-LKKT point of
$D$.

\section{ The main result}

\begin{proposition}
\label{pei}Let $(\overline{x},\overline{\lambda},\overline{\sigma}%
)\in\mathbb{R}^{n}\times\mathbb{R}^{m}\times\operatorname*{int}I^{\ast}$ be a
$J$-LKKT point of $\Xi$ such that $\overline{\sigma}_{k}=0$ for $k\in Q$.

\emph{(i)} Then $\overline{\lambda}\in\Gamma_{J}$, $(\overline{\lambda
},\overline{\sigma})\in T_{Q,\operatorname{col}}$, $\left\langle
\overline{\lambda},\nabla_{\lambda}\Xi(\overline{x},\overline{\lambda
},\overline{\sigma})\right\rangle =0$, $L(\overline{x},\overline{\lambda}%
)=\Xi(\overline{x},\overline{\lambda},\overline{\sigma})$, $\nabla
_{x}L(\overline{x},\overline{\lambda})=0$, and (\ref{r-px2pd}) holds.

\emph{(ii)} Moreover, assume that $Q_{0}^{c}\subset M_{\neq}(\overline
{\lambda})$. Then $\nabla_{\lambda}L(\overline{x},\overline{\lambda}%
)=\nabla_{\lambda}\Xi(\overline{x},\overline{\lambda},\overline{\sigma})$,
$(\overline{x},\overline{\lambda})$ is a $J$-LKKT point of $L$ and
$\overline{x}\in X_{J\cup Q_{0}^{c}}$.

\emph{(iii)} Furthermore, assume that $\overline{\lambda}_{j}>0$ for all $j\in
Q_{0}^{c}$ and $G(\overline{\lambda},\overline{\sigma})\succeq0$. Then
$\overline{x}\in X_{J\cup Q_{0}^{c}}\subset X_{J}\subset X_{J\cap Q}$,
$(\overline{\lambda},\overline{\sigma})\in T_{Q,\operatorname{col}}^{J+}$,
and
\begin{equation}
f(\overline{x})=\inf_{x\in X_{J\cap Q}}f(x)=\Xi(\overline{x},\overline
{\lambda},\overline{\sigma})=L(\overline{x},\overline{\lambda})=\sup
_{(\lambda,\sigma)\in T_{Q,\operatorname{col}}^{J+}}D(\lambda,\sigma
)=D(\overline{\lambda},\overline{\sigma}); \label{r-minmaxei}%
\end{equation}
moreover, if $G(\overline{\lambda},\overline{\sigma})\succ0$ then
$\overline{x}$ is the unique global solution of problem $(P_{J\cap Q})$.
\end{proposition}

Proof. (i) Because $(\overline{x},\overline{\lambda},\overline{\sigma})$ is a
$J$-LKKT point, from its very definition we have that $\overline{\lambda}%
\in\Gamma_{J}$, $\left\langle \overline{\lambda},\nabla_{\lambda}\Xi
(\overline{x},\overline{\lambda},\overline{\sigma})\right\rangle =0$,
$\nabla_{x}\Xi(\overline{x},\overline{\lambda},\overline{\sigma})=0$ and
$\nabla_{\sigma}\Xi(\overline{x},\overline{\lambda},\overline{\sigma})=0$.
Using Lemma \ref{lem-nXiL} and we obtain that $\nabla_{x}L(\overline
{x},\overline{\lambda})=\nabla_{x}\Xi(\overline{x},\overline{\lambda
},\overline{\sigma})=0$ and $L(\overline{x},\overline{\lambda})=\Xi
(\overline{x},\overline{\lambda},\overline{\sigma})$, while using Proposition
\ref{p-perfdual} we get $(\overline{\lambda},\overline{\sigma})\in
T_{Q,\operatorname{col}}$ and that (\ref{r-px2pd}) holds.

(ii) Because $Q_{0}^{c}\subset M_{\neq}(\overline{\lambda})$ we get
$\nabla_{\lambda}L(\overline{x},\overline{\lambda})=\nabla_{\lambda}%
\Xi(\overline{x},\overline{\lambda},\overline{\sigma})$ by Lemma
\ref{lem-nXiL}, and so $(\overline{x},\overline{\lambda})$ is a $J$-LKKT point
of $L$ because $(\overline{x},\overline{\lambda},\overline{\sigma})$ is a
$J$-LKKT point of $\Xi$. Hence $g_{j}(\overline{x})=0$ for $j\in J$, and
$\overline{\lambda}_{j}g_{j}(\overline{x})=0$, $g_{j}(\overline{x})\leq0$ for
$j\in J^{c}$. Taking into account that $Q_{0}^{c}\subset M_{\neq}%
(\overline{\lambda})$, the preceding condition shows that $g_{j}(\overline
{x})=0$ for $j\in Q_{0}^{c,}$ and so $\overline{x}\in X_{J\cup Q_{0}^{c}}$.

(iii) Our hypothesis shows that $Q_{0}^{c}\subset M_{\neq}(\overline{\lambda
})$. From (i) and (ii) we have that $\overline{\lambda}\in\Gamma_{J}$,
$(\overline{\lambda},\overline{\sigma})\in T_{Q,\operatorname{col}}$,
$\overline{x}\in X_{J\cup Q_{0}^{c}}\subset X_{J}\subset X_{J\cap Q}$;
moreover, $\overline{\lambda}\in\Gamma_{J\cap Q}$ because $\overline{\lambda
}_{j}\geq0$ for $j\in J^{c}\cup Q_{0}^{c}=(J\cap Q)^{c}$, and so
$(\overline{\lambda},\overline{\sigma})\in T_{Q,\operatorname{col}}^{J+}$.
Using now Lemma \ref{lem-xiL}, obvious inequalities, (\ref{r-pd2}), and (i),
as well as the obvious inclusion $T_{J,Q\operatorname{col}}^{+}\subset
\Gamma_{J\cap Q}\times I_{J,Q}$ with $I_{J,Q}$ defined in (\ref{r-ikss}), we
get
\begin{align*}
f(\overline{x})  &  \geq\inf_{x\in X_{J\cap Q}}f(x)=\inf_{x\in X_{J\cap Q}%
}\sup_{\lambda\in\Gamma_{J\cap Q}}L(x,\lambda)=\inf_{x\in X_{J\cap Q}}%
\sup_{(\lambda,\sigma)\in\Gamma_{J\cap Q}\times I_{J,Q}}\Xi(x,\lambda
,\sigma)\\
&  \geq\inf_{x\in X_{J\cap Q}}\sup_{(\lambda,\sigma)\in
T_{J,Q\operatorname{col}}^{+}}\Xi(x,\lambda,\sigma)\geq\sup_{(\lambda
,\sigma)\in T_{J,Q\operatorname{col}}^{+}}\inf_{x\in X_{J\cap Q}}\Xi
(x,\lambda,\sigma)\\
&  \geq\sup_{(\lambda,\sigma)\in T_{Q\operatorname{col}}^{J+}}\inf
_{x\in\mathbb{R}^{n}}\Xi(x,\lambda,\sigma)=\sup_{(\lambda,\sigma)\in
T_{Q\operatorname{col}}^{J+}}D(\lambda,\sigma)\geq D(\overline{\lambda
},\overline{\sigma}),
\end{align*}
which implies (\ref{r-minmaxei}) by (i).

Assume, moreover, that $G(\overline{\lambda},\overline{\sigma})\succ0$; hence
$(\overline{\lambda},\overline{\sigma})\in T_{Q}^{J+}$. Consider $x\in
X_{J\cap Q}\setminus\{\overline{x}\}$. Using the strict convexity of
$\Xi(\cdot,\overline{\lambda},\overline{\sigma})$ and Lemma \ref{lem-xiL} we
get $f(\overline{x})=\Xi(\overline{x},\overline{\lambda},\overline{\sigma
})<\Xi(x,\overline{\lambda},\overline{\sigma})\leq L(x,\overline{\lambda})\leq
f(x)$. It follows that $\overline{x}$ is the unique global solution of
$(P_{J\cap Q})$ [and $(P_{J})$, too]. \hfill$\square$

\medskip

The variant of Proposition \ref{pei} in which $Q$ is not taken into
consideration, that is the case when one does not observe that $V_{k}%
\circ\Lambda_{k}=0$ for some $k$, is much weaker; however, the conclusions
coincide for $Q=\{0\}$.

\begin{proposition}
\label{pei-i}Let $(\overline{x},\overline{\lambda},\overline{\sigma}%
)\in\mathbb{R}^{n}\times\mathbb{R}^{m}\times\operatorname*{int}I$ be a
$J$-LKKT point of $\Xi$.

\emph{(i)} Then $\overline{\lambda}\in\Gamma_{J}$, $(\overline{\lambda
},\overline{\sigma})\in T_{\operatorname{col}}$, $\left\langle \overline
{\lambda},\nabla_{\lambda}\Xi(\overline{x},\overline{\lambda},\overline
{\sigma})\right\rangle =0$, $L(\overline{x},\overline{\lambda})=\Xi
(\overline{x},\overline{\lambda},\overline{\sigma})$, $\nabla_{x}%
L(\overline{x},\overline{\lambda})=0$, and (\ref{r-px2pd}) holds.

\emph{(ii)} Assume that $M_{\neq}(\overline{\lambda})=\overline{1,m}$. Then
$\nabla_{\lambda}L(\overline{x},\overline{\lambda})=\nabla_{\lambda}%
\Xi(\overline{x},\overline{\lambda},\overline{\sigma})=0$, whence
$(\overline{x},\overline{\lambda},\overline{\sigma})$ is a critical point of
$\Xi$, $(\overline{x},\overline{\lambda})$ is a critical point of $L$, and
$\overline{x}\in X_{e}\subset X_{J}\subset X_{i}$.

\emph{(iii)} Assume that $\overline{\lambda}\in\mathbb{R}_{++}^{m}$ and
$G(\overline{\lambda},\overline{\sigma})\succeq0$. Then $\overline{x}\in
X_{e}$, $(\overline{\lambda},\overline{\sigma})\in T_{\operatorname{col}}^{+}$
and
\[
f(\overline{x})=\inf_{x\in X_{i}}f(x)=\Xi(\overline{x},\overline{\lambda
},\overline{\sigma})=L(\overline{x},\overline{\lambda})=\sup_{(\lambda
,\sigma)\in T_{\operatorname{col}}^{+}}D(\lambda,\sigma)=D(\overline{\lambda
},\overline{\sigma});
\]
moreover, if $G(\overline{\lambda},\overline{\sigma})\succ0$ then
$(\overline{\lambda},\overline{\sigma})\in T^{+}$ and $\overline{x}$ is the
unique global solution of problem $(P_{i})$.
\end{proposition}

The remark below refers to the case $Q=\emptyset$. A similar remark (but a bit
less dramatic) is valid for $Q_{0}\neq\emptyset$.

\begin{remark}
\label{rem-cdt}It is worth observing that given the functions $f$, $g_{1}$,
..., $g_{m}$ of type $q+V\circ\Lambda$ with $q,\Lambda$ quadratic functions
and $V\in\Gamma_{sc}$, for any choice of $J\subset\overline{1,m}$ one finds
the same $\overline{x}$ using Proposition \ref{pei-i}~(iii). So, in practice,
if one wishes to solve one of the problems $(P_{e})$, $(P_{i})$ or $(P_{J})$
using CDT, it is sufficient to find those critical points $(\overline
{x},\overline{\lambda},\overline{\sigma})$ of $\Xi$ such that $\overline
{\lambda}\in\mathbb{R}_{++}^{m}$ and $G(\overline{\lambda},\overline{\sigma
})\succ0$; if we are successful, $\overline{x}\in X_{e}$ and $\overline{x}$ is
the unique solution of $(P_{i})$, and so $\overline{x}$ is also solution for
all problems $(P_{J})$ with $J\subset\overline{1,m}$; moreover, $(\overline
{\lambda},\overline{\sigma})$ is a global maximizer of $D$ on
$T_{\operatorname{col}}^{+}$.
\end{remark}

The next example shows that the condition $Q_{0}^{c}\subset M_{\neq}%
(\overline{\lambda})$ is essential for $\overline{x}$ to be a feasible
solution of problem $(P_{J})$; moreover, it shows that, unlike the quadratic
case (see \cite[Prop.~9]{Zal:18b}), it is not possible to replace
$T_{Q\operatorname{col}}^{J+}$ by $\{(\lambda,\sigma)\in T_{\operatorname{col}%
}\mid\lambda\in\Gamma_{J},~G(\lambda,\sigma)\succeq0\}$ in (\ref{r-minmaxei}).
The problem is a particular case of the one considered in \cite[Ex.\ 1]%
{LatGao:16}, \textquotedblleft which is very simple, but important in both
theoretical study and real-world applications since the constraint is a
so-called double-well function, the most commonly used nonconvex potential in
physics and engineering sciences [7]\textquotedblright;\footnote{The reference
\textquotedblleft\lbrack7]\textquotedblright\ is \textquotedblleft Gao, D.Y.:
Nonconvex semi-linear problems and canonical duality solutions, in Advances in
Mechanics and Mathematics II. In: Gao, D.Y., Ogden R.W. (eds.), pp. 261--311.
Kluwer Academic Publishers (2003)".} more precisely, $q:=1$, $c:=6$, $d:=4$,
$e:=2$.

\begin{example}
\label{ex1}Let us take $n=m=1$, $J\subset\{1\}$, $q_{0}(x):=\frac{1}{2}%
x^{2}-6x$, $\Lambda_{1}(x):=\frac{1}{2}x^{2}-4$, $q_{1}(x):=\Lambda_{0}%
(x):=0$, $V_{0}(t):=V_{1}(t)+2:=\tfrac{1}{2}t^{2}$ for $x,t\in\mathbb{R}$.
Then $f(x)=\frac{1}{2}x^{2}-6x$ and $g_{1}(x)=\tfrac{1}{2}\left(  \frac{1}%
{2}x^{2}-4\right)  ^{2}-2$. Hence $Q=\{0\}$ (whence $Q_{0}=\emptyset$) and
$X_{e}=\{-2\sqrt{3},2\sqrt{3},-2,2\}\subset\lbrack-2\sqrt{3},-2]\cup
\lbrack2,2\sqrt{3}]=X_{i}$.
\[
\Xi(x;\lambda;\sigma_{0},\sigma)=\tfrac{1}{2}x^{2}-6x-\tfrac{1}{2}\sigma
_{0}^{2}+\lambda\left[  \sigma_{1}\left(  \tfrac{1}{2}x^{2}-4\right)
-\tfrac{1}{2}\sigma_{1}^{2}-2\right]  .
\]
We have that $G(\lambda,\sigma)=1+\lambda\sigma_{1}$, $T_{\operatorname{col}%
}=T=\{(\lambda,\sigma)\in\mathbb{R}\times\mathbb{R}^{2}\mid1+\lambda\sigma
_{1}\neq0\}$ and
\[
D(\lambda;\sigma_{0},\sigma_{1})=-\frac{18}{1+\lambda\sigma_{1}}-\tfrac{1}%
{2}\sigma_{0}^{2}-\lambda\left(  \tfrac{1}{2}\sigma_{1}^{2}+4\sigma
_{1}+2\right)  .
\]
The critical points of $\Xi$ are $\left(  2;-1;(0,-2)\right)  $, $\left(
-2;2;(0,-2)\right)  $, $\left(  6;0;(0,14+8\sqrt{3})\right)  $,
$\big(6;0;(0,14+8\sqrt{3})\big)$, $\left(  -2\sqrt{3};-\frac{1}{2}\sqrt
{3}-\frac{1}{2};(0,2)\right)  $, $\left(  2\sqrt{3};\frac{1}{2}\sqrt{3}%
-\frac{1}{2};(0,2)\right)  $, and so $1+\overline{\lambda}\overline{\sigma
}_{1}\in\{3,-3,1,-\sqrt{3},\sqrt{3}\}$ for $(\overline{x},\overline{\lambda
},\overline{\sigma})$ critical point of $\Xi$, whence $(\overline{\lambda
},\overline{\sigma})$ $(\in T)$ is critical point of $D$ by Lemma \ref{lem1}.
For $\overline{\lambda}=0$ the corresponding $\overline{x}$ $(=6)$ is not in
$X_{i}\supset X_{e}$; in particular, $(\overline{x},\overline{\lambda})$ is
not a critical point of $L$. For $\overline{\lambda}\neq0$, Proposition
\ref{pei} says that $(\overline{x},\overline{\lambda})$ is a critical point of
$L$; in particular $\overline{x}\in X_{e}$. For $\overline{\lambda}%
\in\{2,-\frac{1}{2}\sqrt{3}-\frac{1}{2}\}$, $1+\overline{\lambda}%
\overline{\sigma}_{1}<0$, and so Proposition \ref{pei} says nothing about the
optimality of $\overline{x}$ or $(\overline{\lambda},\overline{\sigma})$; in
fact, for $\overline{\lambda}=-\frac{1}{2}\sqrt{3}-\frac{1}{2}$, the
corresponding $\overline{x}$ $(=-2\sqrt{3})$ is the global maximizer of $f$ on
$X_{e}$. For $\overline{\lambda}:=\frac{1}{2}\sqrt{3}-\frac{1}{2}>0$,
$1+\overline{\lambda}\overline{\sigma}_{1}=\sqrt{3}>0$, and so Proposition
\ref{pei} says that $\overline{x}=2\sqrt{3}$ $(\in X_{e})$ is the global
solution of $(P_{i})$, and $(\overline{\lambda},\overline{\sigma})=\left(
\frac{1}{2}\sqrt{3}-\frac{1}{2};(0,2)\right)  $ is the global maximizer of $D$
on $T_{\operatorname{col}}^{+}=T^{+}=\{(\lambda,\sigma)\in\mathbb{R}_{+}%
\times\mathbb{R}^{2}\mid1+\lambda\sigma_{1}>0\}$. For $\overline{\lambda}=-1$,
$1+\overline{\lambda}\overline{\sigma}_{1}=3>0$, but $(\overline{\lambda
},\overline{\sigma})$ is not a local extremum of $D$, as easily seen taking
$\sigma_{0}:=0$, $(\lambda,\sigma_{1}):=(t-1,t-2)$ with $\left\vert
t\right\vert $ sufficiently small.
\end{example}

When $Q=\overline{0,m}$ problem $(P_{J})$ reduces to the quadratic problem
with equality and inequality quadratic constraints considered in
\cite[$(P_{J})$]{Zal:18b}, which is denoted here by $(P_{J}^{q})$. Of course,
in this case $X=X_{0}=\mathbb{R}^{n}$, and so
\[
\Xi(x,\lambda,\sigma)=L(x,\lambda)-\tfrac{1}{2}\sum\nolimits_{k=0}^{m}%
\lambda_{k}\sigma_{k}^{2}\quad(x\in\mathbb{R}^{n},~\lambda\in\mathbb{R}%
^{m},~\sigma\in\mathbb{R}\times\mathbb{R}^{m})
\]
with $\lambda_{0}:=1$. It follows that
\begin{gather*}
\nabla_{x}\Xi(x,\lambda,\sigma)=\nabla_{x}L(x,\lambda),\quad\nabla_{\sigma}%
\Xi(x,\lambda,\sigma)=-\left(  \lambda_{k}\sigma_{k}\right)  _{k\in
\overline{0,m}},\\
\nabla_{\lambda}\Xi(x,\lambda,\sigma)=\nabla_{\lambda}L(x,\lambda)-\tfrac
{1}{2}\left(  \sigma_{j}^{2}\right)  _{j\in\overline{1,m}}=\left(
q_{j}(x)-\tfrac{1}{2}\sigma_{j}^{2}\right)  _{j\in\overline{1,m}}.
\end{gather*}

Moreover, $G(\lambda,\sigma)=A(\lambda)$, $F(\lambda,\sigma)=b(\lambda)$,
$E(\lambda,\sigma)=c(\lambda)$, and so $T=Y\times\mathbb{R}^{1+m}$,
$T_{\operatorname{col}}=Y_{\operatorname{col}}\times\mathbb{R}^{1+m}$,
$D(\lambda,\sigma)=D(\lambda)-\tfrac{1}{2}\sum_{k=0}^{m}\lambda_{k}\sigma
_{k}^{2}$, where $A(\lambda)$, $b(\lambda)$, $c(\lambda)$, $Y$,
$Y_{\operatorname{col}}$, $D$ are introduced in \cite{Zal:18b}; we set
$D_{L}:=D$ in the present case. Applying Proposition \ref{pei-i} for this case
we get the next result.

\begin{corollary}
\label{c-pei-i}Let $(\overline{x},\overline{\lambda})\in\mathbb{R}^{n}%
\times\mathbb{R}^{m}$ be a $J$-LKKT point of $L$.

\emph{(i)} Then $\overline{\lambda}\in Y_{\operatorname{col}}^{J}%
:=Y_{\operatorname{col}}\cap\Gamma_{J}$, $\left\langle \overline{\lambda
},\nabla_{\lambda}L(\overline{x},\overline{\lambda})\right\rangle =0$, and
$q_{0}(\overline{x})=L(\overline{x},\overline{\lambda})=D_{L}(\overline
{\lambda})$.

\emph{(ii)} Assume that $M_{\neq}(\overline{\lambda})=\overline{1,m}$. Then
$\nabla_{\lambda}L(\overline{x},\overline{\lambda})=0$, and so $(\overline
{x},\overline{\lambda})$ is a critical point of $L$, and $\overline{x}\in
X_{e}\subset X_{J}\subset X_{i}$.

\emph{(iii)} Assume that $\overline{\lambda}\in\mathbb{R}_{++}^{m}$ and
$A(\overline{\lambda})\succeq0$. Then $\overline{x}\in X_{e}$, $\overline
{\lambda}\in Y_{\operatorname{col}}^{+}$ and%
\[
q_{0}(\overline{x})=\inf_{x\in X_{i}}q_{0}(x)=L(\overline{x},\overline
{\lambda})=\sup_{\lambda\in Y_{\operatorname{col}}^{i+}}D_{L}(\lambda
)=D_{L}(\overline{\lambda});
\]
moreover, if $A(\overline{\lambda})\succ0$ then $\overline{\lambda}\in Y^{i+}$
and $\overline{x}$ is the unique global solution of problem $(P_{i})$.
\end{corollary}

However, applying Proposition \ref{pei} we get assertion (i) and last part of
assertion (ii) of \cite[Prop.~9]{Zal:18b}.

As seen in \cite[Prop.~9]{Zal:18b} the most part of the results obtained by
DY\ Gao and his collaborators for quadratic minimization problems are very far
from those obtained studying directly those quadratic problems. In this sense
it is worth quoting the following remark from the very recent Ruan and Gao's
paper \cite{RuaGao:17b}:

\begin{quote}
``\emph{Remark 1}. As we have demonstrated that by the generalized canonical
duality (32), all KKT conditions can be recovered for both equality and
inequality constraints. Generally speaking, the nonzero Lagrange multiplier
condition for the linear equality constraint is usually ignored in
optimization textbooks. But it can not be ignored for nonlinear constraints.
It is proved recently [26] that the popular augmented Lagrange multiplier
method can be used mainly for linear constrained problems. Since the
inequality constraint $\mu\not =0$ produces a nonconvex feasible set
$\mathcal{E}_{a}^{\ast}$, this constraint can be replaced by either $\mu<0$ or
$\mu>0$. But the condition $\mu<0$ is corresponding to $y\circ(y-e_{K})\geq0$,
this leads to a nonconvex open feasible set for the primal problem. By the
fact that the integer constraints $y_{i}(y_{i}-1)=0$ are actually a special
case (boundary) of the boxed constraints $0\leq y_{i}\leq1$, which is
corresponding to $y\circ(y-e_{K})\geq0$, we should have $\mu>0$ (see [8] and
[12, 16]). In this case, the KKT condition (43) should be replaced by

$\mu>0,~~y\circ(y-e_{K})\leq0,~~\mu^{T}[y\circ(y-e_{K})]=0.\quad$ (47)

\noindent Therefore, as long as $\mu\neq0$ is satisfied, the complementarity
condition in (47) leads to the integer condition $y\circ(y-e_{K})=0$.
Similarly, the inequality $\tau\neq0$ can be replaced by $\tau>0$%
."\footnote{The reference \textquotedblleft\lbrack26]" mentioned in
\cite[Rem.~1]{RuaGao:17b} is the item \cite{LatGao:16} from our bibliography,
the others being the following: \textquotedblleft8.~Fang, S.C., Gao, D.Y.,
Sheu, R.L., Wu, S.Y.: Canonical dual approach to solving 0--1 quadratic
programming problems.\ J. Ind.\ Manag.\ Optim.\ 4(4), 125--142 (2008)",
\textquotedblleft12.~Gao, D.Y.: Solutions and optimality criteria to box
constrained nonconvex minimization problem.\ J.\ Ind.\ Manag.\ Optim.\ 3(2),
293--304 (2007)", and \textquotedblleft16.~Gao, D.Y., Ruan, N.: Solutions to
quadratic minimization problems with box and integer constraints. J.
Glob.\ Optim.\ 47, 463--484 (2010)", respectively.}
\end{quote}

In fact the positivity of the Lagrange multipliers $\lambda_{j}$ is needed for
recovering the Lagrangian $L$ from $\Xi$ [see (\ref{r-lgk})], while the non
vanishing condition on $\overline{\lambda}_{j}$ is needed to get
$L(\overline{x},\overline{\lambda})=\Xi(\overline{x},\overline{\lambda
},\overline{\sigma})$ and $\nabla_{x}L(\overline{x},\overline{\lambda}%
)=\nabla_{x}\Xi(\overline{x},\overline{\lambda},\overline{\sigma})$ when
$\nabla_{\sigma}\Xi(\overline{x},\overline{\lambda},\overline{\sigma})=0$, as
seen in Lemma \ref{lem-nXiL}. Of course, such conditions are not needed in
quadratic minimization problems, as observed after Corollary \ref{c-pei-i}.

\section{Relations with previous results}

In this section we analyze results obtained by DY Gao and his collaborators in
papers dedicated to constrained optimization problems. Because the quadratic
problems (with quadratic constraints) are discussed in \cite{Zal:18b}, we
discuss only those constrained optimization problems with non quadratic
objective function or with at least one non quadratic constraint. In the
survey paper \cite{GaoRuaLat:17} (the same as \cite{GaoRuaLat:16}) there are
mentioned the following papers: \cite{GaoYan:08}, \cite{GaoRuaShe:09},
\cite{LatGao:16}, \cite{MorGao:17} (with its preprint version \cite{MorGao:12}%
); besides these papers we add the retracted version \cite{MorGao:16} of
\cite{MorGao:17}, and \cite{RuaGao:17}.

A detailed discussion of \cite{GaoYan:08} was done in \cite{VoiZal:11}; we
discuss the corrected versions \cite{MorGao:12}--\cite{MorGao:17} of
\cite{GaoYan:08} at the end of this section.

\medskip The problem considered by Gao, Ruan and Sherali in
\cite{GaoRuaShe:09} is of type $(P_{i})$, that is $J=\emptyset$ with our
notation, with $f$ a quadratic function. Taking $q_{j}:=0$ for $j\in
\overline{1,m}$, our problem $(P_{i})$ is a particular case of the problem
$(\mathcal{P})$ from \cite{GaoRuaShe:09}. In this framework, that is
$V_{0}(y)=\tfrac{1}{2}y^{2}$, $V_{j}\in\Gamma_{sc}$, $\Lambda_{0}:=q_{j}:=0$
for $j\in\overline{1,m}$, with our notations, we mention only the following
result of \cite{GaoRuaShe:09}.

\smallskip\textquotedblleft Theorem 2 (Global Optimality Condition)".
\emph{Let $(\overline{x},\overline{\lambda},\overline{\sigma})\in
X\times\mathbb{R}_{+}^{m}\times I^{\ast}$ be a critical point of $\Xi$. If
$G(\overline{\lambda},\overline{\sigma})\succeq0$, then $(\overline{\lambda
},\overline{\sigma})$ is a global maximizer of $D$ on $T_{\operatorname{col}%
}^{+}$, $\overline{x}$ is a global minimizer of $f$ on $X_{i}$ and
$f(\overline{x})=\min_{x\in X_{i}}f(x)=\max_{(\lambda,\sigma)\in
T_{\operatorname{col}}^{+}}D(\lambda,\sigma)=D(\overline{\lambda}%
,\overline{\sigma}).$}

\smallskip This theorem is false because in the mentioned conditions
$\overline{x}$ is not necessarily in $X_{i}$, as Example \ref{ex1} shows.
Indeed, $\left(  6;0;(0,14+8\sqrt{3})\right)  $ is a critical point of $\Xi$,
but $6\notin X_{i}$. It follows that also \textquotedblleft Theorem 1
(Complementary-Dual Principle)" and \textquotedblleft Theorem 3 (Triality
Theory)" of \cite{GaoRuaShe:09} are false because $(\overline{\lambda
},\overline{\sigma})=\left(  0;(0,14+8\sqrt{3})\right)  $ is a critical point
of $D$ (by Lemma \ref{lem1}), but the assertion \textquotedblleft$\overline
{x}$ is a KKT point of $(\mathcal{P})$" is not true.

It is shown in \cite[Ex.~6]{VoiZal:11b} that the \textquotedblleft double-min
or double-max" duality of \cite[Theorem 3 (Triality Theory)]{GaoRuaShe:09},
that is its assertion in the case $G(\overline{\lambda},\overline{\sigma
})\prec0$, is also false.

\medskip

The problem considered by Latorre and Gao in \cite{LatGao:16} is of type
$(P_{J})$ in which $\Lambda_{k}$ are quadratic and $V_{k}$ are
\textquotedblleft differentiable canonical functions". In our framework
(which, apparently, is more restrictive) and with our notations, the following
set is used in \cite{LatGao:16}:%
\[
\mathcal{S}_{0}:=\{\lambda\in\mathbb{R}^{m}\mid\left[  \forall j\in
J:\lambda_{j}\neq0\right]  ~\wedge~\left[  \forall j\in J^{c}:\lambda_{j}%
\geq0\right]  \}\subset\Gamma_{J}.
\]
The motivation for defining $\mathcal{S}_{0}$ like this is given in the
following text from \cite[p.~1767]{LatGao:16}: \textquotedblleft From the
second and third equation in the (10), it is clear that in order to enforce
the constrain $h(x)=0$, the dual variables $\mu_{i}$ must be not zero for
$i=1,...,p$. \emph{This is a special complementarity condition for equality
constrains, generally not mentioned in many textbooks.} However, the implicit
constraint $\mu\neq0$ is important in nonconvex optimization. Let $\sigma
_{0}=(\lambda,\mu)$. The dual feasible spaces should be defined as
$\mathcal{S}_{0}$ ...".\footnote{The emphasized text can be found also in
\cite[p.~NP26]{GaoRuaLat:16} and \cite[p.~33]{GaoRuaLat:17}. One must also
observe that for Latorre and Gao $\mu\neq0$ is equivalent to \textquotedblleft%
$\mu_{j}\neq0$ $\forall j=1,...,p$", and $(\lambda,\mu)\in\mathbb{R}^{m\times
p}$ if $\lambda\in\mathbb{R}^{m}$ and $\mu\in\mathbb{R}^{p}.$}

Besides the set $\mathcal{S}_{0}$ mentioned above, the following sets are also
considered in \cite{LatGao:16}:
\[
\mathcal{S}_{1}:={\textstyle\prod\nolimits_{k=0}^{m}}\operatorname*{dom}%
V_{k}=I^{\ast},\quad\mathcal{S}_{a}:=T_{\operatorname{col}}\cap\left(
\mathcal{S}_{0}\times\mathcal{S}_{1}\right)  ,\quad\mathcal{S}_{a}%
^{+}:=\left\{  (\lambda,\sigma)\in T^{+}\mid J\subset M_{\neq}(\lambda
)\right\}  .
\]

In this context the main results of \cite{LatGao:16} are the following.

\smallskip\textquotedblleft Theorem 1 (Complementarity Dual Principle)".
\emph{Let $(\overline{x},\overline{\lambda},\overline{\sigma})\in
X\times\mathbb{R}_{+}^{m}\times I^{\ast}$ be a critical point of $\Xi$. Then
$\overline{x}$ is a $J$-KKT of $(P_{J})$, $(\overline{\lambda},\overline
{\sigma})$ is a $J$-LKKT point of $D$ and $f(\overline{x})=\Xi(\overline
{x},\overline{\lambda},\overline{\sigma})=D(\overline{\lambda},\overline
{\sigma})$.}

\smallskip\textquotedblleft Theorem 2 (Global Optimality Conditions)".
\emph{Let $(\overline{x},\overline{\lambda},\overline{\sigma})\in
X\times\mathbb{R}_{+}^{m}\times I^{\ast}$ be a critical point of $\Xi$ with
$(\overline{\lambda},\overline{\sigma})\in\mathcal{S}_{a}^{+}$. If
$\mathcal{S}_{a}^{+}$ is convex then $(\overline{\lambda},\overline{\sigma})$
is the global maximizer of $D$ on $\mathcal{S}_{a}^{+}$ and $\overline{x}$ is
the global minimizer of $f$ on $X_{J}$, that is $f(\overline{x})=\min_{x\in
X_{J}}f(x)=\max_{(\lambda,\sigma)\in\mathcal{S}_{a}^{+}}D(\lambda
,\sigma)=D(\overline{\lambda},\overline{\sigma})$.}

\smallskip Note first that it is not clear what is meant by $J$-LKKT point of
$D$ (called KKT point) in \cite[Th.~1]{LatGao:16} when $(\overline{\lambda
},\overline{\sigma})\notin T$. As in the case of \cite[Th.~1]{GaoRuaShe:09},
Example \ref{ex1} shows that \cite[Th.~1]{LatGao:16} is false because
$\big(6;0;(0,14+8\sqrt{3})\big)$ is a critical point of $\Xi$, but $6\notin
X_{J}$ $(=X_{i})$; even without assuming that $\mathcal{S}_{a}^{+}$ is convex
in \cite[Th.~2]{LatGao:16}, for the same reason, this theorem is false.

Having in view that there are not nonempty open convex subsets $C\subset
\mathbb{R}^{2}$ such that the mapping $C\ni(u,v)\mapsto uv\in\mathbb{R}$ is
convex, the hypothesis that $\mathcal{S}_{a}^{+}$ is convex in the statement
of \cite[Th.~2]{LatGao:16} is very strong. Moreover, it is not clear how this
hypothesis is used in the proof of \cite[Th.~2]{LatGao:16}.\footnote{It is
worth quoting DY Gao's comment from \cite[p.~19]{Gao:16} on our remark from
\cite[p.~1783]{Zal:16} that the proof of \cite[Th.~2]{LatGao:16} is not
convincing: \textquotedblleft Regarding the so-called \textquotedblleft not
convincing proof\textquotedblright, serious researcher should provide either a
convincing proof or a disproof, rather than a complaint. Note that the
canonical dual variables $\sigma_{0}$ and $\sigma_{1}$ are in two different
levers (scales) with totally different physical units$^{14}$, it is completely
wrong to consider $(\sigma_{0},\sigma_{1})$ as one vector and to discuss the
concavity of $\Xi_{1}\left(  x,(\cdot,\cdot)\right)  $ on $\mathcal{S}_{a}%
^{+}$. The condition \textquotedblleft$\mathcal{S}_{a}^{+}$ is
convex\textquotedblright\ in Theorem 2 [5] should be understood in the way
that $\mathcal{S}_{a}^{+}$ is convex in $\sigma_{0}$ and $\sigma_{1}$,
respectively, as emphasized in Remark 1 [5]. Thus, the proof of Theorem 2
given in [5] is indeed convincing by simply using the classical saddle min-max
duality for $(x,\sigma_{0})$ and $(x,\sigma_{1})$, respectively."
\par
Note 14 from the text above is \textquotedblleft Let us consider Example 1 in
[5]. If the unit for $x$ is the meter $(m)$ and for $q$ is $Kg/m$, then the
units for the Lagrange multiplier $\mu$ (dual to the constraint $g(x)=\tfrac
{1}{2}(\tfrac{1}{2}x^{2}-d)^{2}-e$) should be $Kg/m^{3}$ and for $\sigma$
(canonical dual to $\Lambda(x)=\tfrac{1}{2}x^{2}$) should be $Kg/m$,
respectively, so that each terms in $\Xi_{1}(x,\mu,\sigma)$ make physical
sense"; \textquotedblleft\lbrack5]" is our reference \cite{LatGao:16}.}

\medskip

The results established by Ruan and Gao in Sections 3 of \cite{RuaGao:16} and
\cite{RuaGao:17} (which are practically the same) refer to $(P_{i})$ in which
$q_{k}=0$, $\Lambda_{k}$ are G\^{a}teaux differentiable on their domains and
$V_{k}$ are \textquotedblleft canonical functions" for $k\in\overline{0,m}$.
In our framework (which is more restrictive) and with our notations, the
following sets are used in \cite{RuaGao:17}:%
\[
\mathcal{S}_{a}:=T_{\operatorname{col}}\cap\left(  \mathbb{R}_{+}^{m}\times
I^{\ast}\right)  ,\quad\mathcal{S}_{a}^{+}:=\left\{  (\lambda,\sigma)\in
T^{+}\mid M_{\neq}(\lambda)=\overline{1,m}\right\}  .
\]

In this context the results of \cite{RuaGao:16} and \cite{RuaGao:17} we are
interested in are the following.

\smallskip

Theorem 3. \emph{Let $(\overline{x},\overline{\lambda},\overline{\sigma})\in
X\times\mathbb{R}_{+}^{m}\times I^{\ast}$ be a KKT point of $\Xi$. Then
$\overline{x}$ is a KKT of $(P_{i})$, $(\overline{\lambda},\overline{\sigma})$
is a KKT point of $D$ and $f(\overline{x})=\Xi(\overline{x},\overline{\lambda
},\overline{\sigma})=D(\overline{\lambda},\overline{\sigma})$.}

\smallskip Theorem 4. \emph{Let $(\overline{x},\overline{\lambda}%
,\overline{\sigma})\in X\times\mathbb{R}_{+}^{m}\times I^{\ast}$ be a KKT
point of $\Xi$ with $(\overline{\lambda},\overline{\sigma})\in\mathcal{S}%
_{a}^{+}$. If $\mathcal{S}_{a}^{+}$ is convex then $(\overline{\lambda
},\overline{\sigma})$ is a global maximizer of $D$ on $\mathcal{S}_{a}^{+}$
and $\overline{x}$ is a global minimizer of $f$ on $X_{i}$, that is
$f(\overline{x})=\min_{x\in X_{i}}f(x)=\max_{(\lambda,\sigma)\in
\mathcal{S}_{a}^{+}}D(\lambda,\sigma)=D(\overline{\lambda},\overline{\sigma}%
)$.}

\smallskip As in \cite[Th.~1]{LatGao:16}, it is not clear what is meant by KKT
point of $D$ in \cite[Th.~3]{RuaGao:17} when $(\overline{\lambda}%
,\overline{\sigma})\notin T$. As in the case of \cite[Th.~1]{GaoRuaShe:09},
Example \ref{ex1} shows that \cite[Th.~3]{RuaGao:17} is false because
$\big(6;0;(0,14+8\sqrt{3})\big)$ is a critical point of $\Xi$, hence a KKT
point of $\Xi$, but $6\notin X_{i}$. In what concerns \cite[Th.~4]{RuaGao:17},
because $M_{\neq}(\lambda)=\overline{1,m}$, $(\overline{x},\overline{\lambda
},\overline{\sigma})$ is a critical point of $\Xi$ and $\overline{x}\in X_{e}%
$; moreover, in our framework (that is $V_{k}\in\Gamma_{sc}$ for
$k\in\overline{0,m}$), this theorem is true without assuming that
$\mathcal{S}_{a}^{+}$ is convex. Notice that the proof of \cite[Th.~4]%
{RuaGao:17} is not convincing.

\medskip

Morales-Silva and Gao in \cite{MorGao:12}--\cite{MorGao:17} consider the
problem $(\mathcal{P})$ of minimizing $\tfrac{1}{2}\left\Vert y-z\right\Vert
^{2}$ for $x:=(y,z)\in\mathcal{Y}_{c}\times\mathcal{Z}_{c}$ with
$\mathcal{Y}_{c}:=\left\{  y\in\mathbb{R}^{n}\mid h(y)=0\right\}  $ and
$\mathcal{Z}_{c}:=\left\{  z\in\mathbb{R}^{n}\mid h(z)=0\right\}  $, where
$h(y):=\tfrac{1}{2}\left(  \left\langle y,Ay\right\rangle -r^{2}\right)  $ and
$h(z):=\tfrac{1}{2}\alpha\big(\tfrac{1}{2}\left\Vert z-c\right\Vert ^{2}%
-\eta\big)^{2}-\left\langle f,z-c\right\rangle $; here $A\in\mathfrak{S}_{n}$
is positive definite, $c,f\in\mathbb{R}^{n}$ and $\alpha,\eta,r\in(0,\infty)$
are taken such that $h(z)>0$ for every $z\in\mathcal{Z}_{c}$. Of course, this
problem is of type $(P_{e})$ for which Proposition \ref{pei-i} applies.
Because \cite{MorGao:12} is the preprint version of \cite{MorGao:17}, we refer
mostly to \cite{MorGao:17} and \cite{MorGao:16}.\footnote{Excepting
\cite{MorGao:17}, there are very few differences between the papers published
in \cite{GaoLatRua:17} and those having the same title from the retracted
issue of the journal Mathematics and Mechanics of Solids dedicated to CDT.} In
\cite{MorGao:17} one considers the sets

\textquotedblleft$\mathcal{S}_{a}=\{(\lambda,\mu,\varsigma)\in\mathbb{R}%
\times\mathbb{R}\times\mathcal{V}_{a}^{\ast}:(1+\mu\varsigma)(I+\lambda A)-I$
is invertible$\}.\quad\mathcal{(}10\mathcal{)}$

$\mathcal{S}_{a}^{+}=\{(\lambda,\mu,\varsigma)\in\mathcal{S}_{a}:I+\lambda
A\succ0$ and $(1+\mu\varsigma)(I+\lambda A)-I\succ0\}.\quad\mathcal{(}%
19\mathcal{)}$"

\noindent where $\mathcal{V}_{a}^{\ast}:=[-\alpha\eta,\infty)$, and one states
the following results:

\smallskip\textquotedblleft Theorem 1 (Complementary-dual principle). \emph{If
$(\overline{x},\overline{\lambda},\overline{\mu},\overline{\varsigma})$ is a
stationary point of $\Xi$ such that $(\overline{\lambda},\overline{\mu
},\overline{\varsigma})\in\mathcal{S}_{a}$ then $\overline{x}$ is a critical
point of $(\mathcal{P})$ with $\overline{\lambda}$ and $\overline{\mu}$ its
Lagrange multipliers, $(\overline{\lambda},\overline{\mu},\overline{\varsigma
})$ is a stationary point of $\Pi^{d}$ and $\Pi(\overline{x})=L(\overline
{x},\overline{\lambda},\overline{\mu})=\Xi(\overline{x},\overline{\lambda
},\overline{\mu},\overline{\varsigma})=\Pi^{d}(\overline{\lambda}%
,\overline{\mu},\overline{\varsigma}).\quad(17)$}"

\smallskip\textquotedblleft Theorem 2. \emph{Suppose that $(\overline{\lambda
},\overline{\mu},\overline{\varsigma})\in\mathcal{S}_{a}^{+}$ is a stationary
point of $\Pi^{d}$ with $\overline{\mu}\geq0$. Then $\overline{x}$ defined by
(11) is the only global minimizer of $\Pi$ on $\mathcal{X}_{c}$, and
$\Pi(\overline{x})=\min_{x\in\mathcal{X}_{c}}\Pi(x)=\max_{(\lambda
,\mu,\varsigma)\in\mathcal{S}_{a}^{+}}\Pi^{d}(\lambda,\mu,\varsigma)=\Pi
^{d}(\overline{\lambda},\overline{\mu},\overline{\varsigma}).\quad
(20)$\textquotedblright}

\smallskip Theorem 2.2 of \cite{MorGao:12} coincides with \cite[Th.~1]%
{MorGao:17}, while in Theorem 2.3 of \cite{MorGao:12} \textquotedblleft%
$\overline{\mu}\geq0$" and Eq.~(20) from the statement of \cite[Th.~2]%
{MorGao:17} are missing.

\medskip In \cite{MorGao:16}, in the context of the problem $(\mathcal{P})$
above, one considers the sets

\textquotedblleft$\mathcal{S}_{a}=\{(\lambda,\mu,\varsigma)\in\mathbb{R}%
\times\mathbb{R}\times\mathcal{V}_{a}^{\ast}:\lambda\neq0,~\mu\neq
0,~\det\left[  (1+\mu\varsigma)(I+\lambda A)-I\right]  \neq0\}.\quad
\mathcal{(}27\mathcal{)}$

$\mathcal{S}_{c}^{+}=\{(\lambda,\mu,\varsigma)\in\mathcal{S}_{a}:\lambda
>0,\mu>0,~I+\lambda A\succ0$ and $(1+\mu\varsigma)(I+\lambda A)-I\succ
0\}.\quad\mathcal{(}30\mathcal{)}$".

With this new $\mathcal{S}_{a}$, \cite[Th.~2]{MorGao:16} has the same
statement as \cite[Th.~1]{MorGao:17}; moreover, replacing $\mathcal{S}_{a}%
^{+}$ with $\mathcal{S}_{c}^{+}$ in the statement of \cite[Th.~2]{MorGao:17}
one gets the statement of \cite[Th.~3]{MorGao:16}.

\smallskip Notice that there is not a proof of the equality $\max
_{(\lambda,\mu,\varsigma)\in\mathcal{S}_{a}^{+}}\Pi^{d}(\lambda,\mu
,\varsigma)=\Pi^{d}(\overline{\lambda},\overline{\mu},\overline{\varsigma})$
in \cite{MorGao:17}, and there is not a proof of the equality $\max
_{(\lambda,\mu,\varsigma)\in\mathcal{S}_{c}^{+}}\Pi^{d}(\lambda,\mu
,\varsigma)=\Pi^{d}(\overline{\lambda},\overline{\mu},\overline{\varsigma})$
in \cite{MorGao:16}. However, there is a \textquotedblleft proof" of the
equality $\max_{(\lambda,\mu,\varsigma)\in\mathcal{S}_{c}}P^{d}(\lambda
,\mu,\varsigma)=P^{d}(\overline{\lambda},\overline{\mu},\overline{\varsigma})$
from Theorem 2 of \cite{GaoYan:08} even if $\mathcal{S}_{c}$ defined in
\cite[Eq.~(16)]{GaoYan:08} includes $\mathcal{S}_{a}^{+}$ defined in
\cite[Eq.~(19)]{MorGao:17}; see the discussion form \cite[Sect.~2]%
{VoiZal:11}.\footnote{In fact, we did not find a correct proof of the
\textquotedblleft min-max" duality (like Eq.~(20) from \cite{MorGao:17}) in DY
Gao's papers in the case $Q_{0}\neq\emptyset.$}

\medskip

Setting $g_{1}:=h$ and $g_{2}:=g$, we have that $Q=\{0,1\}$ and $J=\{1,2\}$ in
problem $(\mathcal{P})$ of \cite{MorGao:12}--\cite{MorGao:17}. Because
$\mathcal{Y}_{c}\cap\mathcal{Z}_{c}=\emptyset$ and taking into account
\cite[Assertion II, p.~596]{VoiZal:11},\footnote{Referring to \cite[Lem.~1]%
{MorGao:16}, which is a reformulation of \cite[Assertion II]{VoiZal:11}, the
authors say: \textquotedblleft The following Lemma is well known in
mathematical programming (cf.\ Latorre and Gao [12] and Voisei and Zalinescu
[13])"; the references \textquotedblleft\lbrack12]\textquotedblright\ and
\textquotedblleft\lbrack13]" are our items \cite{LatGao:16} and
\cite{VoiZal:11}, respectively. Of course, \cite[Lem.~1]{MorGao:16} is not
\textquotedblleft well known in mathematical programming", being very specific
to the problem considered in \cite{GaoYan:08}. The reference \cite{LatGao:16}
is not mentioned in \cite{MorGao:17} with respect to \cite[Lem.~1]%
{MorGao:16}.} under the hypothesis of \cite[Th.~1]{MorGao:17} one has
$\overline{\lambda}\neq0\neq\overline{\mu}$, and so $M_{\neq}(\overline
{\lambda},\overline{\mu})=\{1,2\}$. Using Proposition \ref{pei-i}~(ii) we
obtain that $(\overline{x},\overline{\lambda},\overline{\mu})$ is a critical
point of $L$ and \cite[Eq.~(17)]{MorGao:17} holds. The conclusion of
\cite[Th.~3]{MorGao:16} is obtained using Proposition \ref{pei-i}~(iii)
[taking into account Corollary \ref{c-lkkt} (ii)]. In what concerns
\cite[Th.~2]{MorGao:17}, its conclusion follows using Proposition
\ref{pei}~(iii) because the condition [$\mathcal{Y}_{c}\cap\mathcal{Z}%
_{c}=\emptyset$ $\wedge$ $\overline{\mu}\geq0$] imply $\overline{\mu}>0$, and
so $Q_{0}^{c}=\{2\}\subset M_{\neq}(\overline{\lambda},\overline{\mu})$.

\medskip Below we show that the equality $\max_{(\lambda,\mu,\varsigma
)\in\mathcal{S}_{a}^{+}}\Pi^{d}(\lambda,\mu,\varsigma)=\Pi^{d}(\overline
{\lambda},\overline{\mu},\overline{\varsigma})$ from \cite[Eq.~(20)]%
{MorGao:17} is not true. For this consider $n:=1$, $A:=1$, $r:=\alpha
:=\eta:=c:=1$ and $f:=\frac{9}{8}\sqrt{2}$; this is a particular case
$(\gamma:=\frac{9}{8}\sqrt{2})$ of the problem $(\mathcal{P})$ considered in
\cite{VoiZal:11}. In this situation (with the calculations and notations from
\cite{VoiZal:11}), the equation $\varsigma^{4}=8\gamma^{2}(\varsigma+1)$ has
the solutions $\overline{\varsigma}:=\varsigma_{1}\in(-1,0)$ (and so
$\overline{\varsigma}+\gamma>0$), and $\varsigma_{2}=3$. Taking $\overline
{\lambda}:=\frac{\overline{\varsigma}^{2}}{2\gamma}>0$ and $\overline{\mu
}:=\frac{\overline{\varsigma}^{2}}{2\gamma^{2}-\overline{\varsigma}^{3}}>0$,
we have that $(\overline{\lambda},\overline{\mu},\overline{\varsigma})$ is a
critical point of $D$ $(=\Pi^{d})$; moreover, $1+\overline{\lambda}>0$ and
$(1+\overline{\lambda})(1+\overline{\mu}\overline{\varsigma})-1=\overline{\mu
}(\gamma+\overline{\varsigma})>0$, and so $(\overline{\lambda},\overline{\mu
},\overline{\varsigma})\in\mathcal{S}_{a}^{+}$, where $\mathcal{S}_{a}^{+}$ is
defined in \cite[Eq.~(19)]{MorGao:17}. In fact, in the present case,
\[
D(\lambda,\mu,\varsigma)=\Pi^{d}(\lambda,\mu,\varsigma)=-\frac{\mu^{2}\left(
\lambda+1\right)  \left(  \varsigma^{3}+2\varsigma^{2}+\gamma^{2}\right)
+\mu\lambda\left(  \varsigma^{2}+\varsigma\lambda+2\varsigma-2\gamma\right)
+\lambda^{2}}{2\left(  \lambda+\varsigma\mu+\varsigma\lambda\mu\right)  }%
\]
for all $(\lambda,\mu,\varsigma)\in\mathcal{S}_{a}$. Applying \cite[Th.~2]%
{MorGao:17} we must have that $\max_{(\lambda,\mu,\varsigma)\in\mathcal{S}%
_{a}^{+}}\Pi^{d}(\lambda,\mu,\varsigma)=\Pi^{d}(\overline{\lambda}%
,\overline{\mu},\overline{\varsigma})$. However, this is not possible because
$\sup_{(\lambda,\mu,\varsigma)\in\mathcal{S}_{a}^{+}}\Pi^{d}(\lambda
,\mu,\varsigma)=\infty$. Indeed, there exists $\tilde{\varsigma}<0$ such that
$\nu:=\tilde{\varsigma}^{3}+2\tilde{\varsigma}^{2}+\gamma^{2}<0$. Then
$(\overline{\lambda},\mu,\tilde{\varsigma})\in\mathcal{S}_{a}^{+}$ for every
$\mu<0$ because $1+\overline{\lambda}>0$ and $(1+\overline{\lambda}%
)(1+\tilde{\varsigma}\mu)-1\geq(1+\overline{\lambda})-1=\overline{\lambda}>0$.
It follows that
\begin{align*}
D(\overline{\lambda},\mu,\tilde{\varsigma})  &  =-\frac{\mu^{2}\left(
\overline{\lambda}+1\right)  \left(  \tilde{\varsigma}^{3}+2\tilde{\varsigma
}^{2}+\gamma^{2}\right)  +\mu\overline{\lambda}\left(  \tilde{\varsigma}%
^{2}+\tilde{\varsigma}\overline{\lambda}+2\tilde{\varsigma}-2\gamma\right)
+\overline{\lambda}^{2}}{2\left(  \overline{\lambda}+\tilde{\varsigma}%
\mu+\tilde{\varsigma}\overline{\lambda}\mu\right)  }\\
&  =-\frac{\mu^{2}\left(  \overline{\lambda}+1\right)  \nu+\mu\overline
{\lambda}\left(  \tilde{\varsigma}^{2}+\tilde{\varsigma}\overline{\lambda
}+2\tilde{\varsigma}-2\gamma\right)  +\overline{\lambda}^{2}}{2\left(
\overline{\lambda}+\tilde{\varsigma}\mu+\tilde{\varsigma}\overline{\lambda}%
\mu\right)  }\quad\forall\mu<0,
\end{align*}
and so
\begin{align*}
\sup_{(\lambda,\mu,\varsigma)\in\mathcal{S}_{a}^{+}}\Pi^{d}(\lambda
,\mu,\varsigma)  &  \geq-\lim_{\mu\rightarrow-\infty}\frac{\mu^{2}\left(
\overline{\lambda}+1\right)  \nu+\mu\overline{\lambda}\left(  \tilde
{\varsigma}^{2}+\tilde{\varsigma}\overline{\lambda}+2\tilde{\varsigma}%
-2\gamma\right)  +\overline{\lambda}^{2}}{2\left(  \left[  \overline{\lambda
}+\left(  \overline{\lambda}+1\right)  \tilde{\varsigma}\mu\right]  \right)
}\\
&  =-\lim_{\mu\rightarrow-\infty}\frac{\mu^{2}\left(  \overline{\lambda
}+1\right)  \nu}{\left(  \overline{\lambda}+1\right)  \tilde{\varsigma}\mu
}=-\lim_{\mu\rightarrow-\infty}\frac{\nu}{\tilde{\varsigma}}\mu=\infty.
\end{align*}

In \cite{VoiZal:11} we provided an example with $n=2$ for which the
solution(s) of problem $(\mathcal{P})$ from \cite{GaoYan:08} (which clearly
always exists) can not be obtained (found) using \cite[Th.~2]{GaoYan:08}; we
concluded that \textquotedblleft the consideration of the function $\Xi$ is
useless, at least for the problem studied in [3]".\footnote{The reference
\textquotedblleft\lbrack3]" is the item \cite{GaoYan:08} from our
bibliography.} In \cite{MorGao:12}--\cite{MorGao:17} the authors sustain that
this is caused by the non uniqueness of the solution of problem $(\mathcal{P}%
)$ from our example, but a solution can be obtained, even in such a case, by
perturbation: \textquotedblleft The combination of the perturbation and the
canonical duality theory is an important method for solving nonconvex
optimization problems which have more than one global optimal solution (see
also [15])."\footnote{The text is quoted from \cite[p.~370]{MorGao:17}; here
the reference \textquotedblleft\lbrack15]" is \textquotedblleft Wu, C., Gao,
D.Y.: Canonical primal-dual method for solving nonconvex minimization
problems. In: Gao, D.Y., Latorre, V., Ruan, N. (eds.) Advances in Canonical
Duality Theory. Springer, Berlin". Note that the same text can be found in
\cite[p.~9]{MorGao:12} without any reference, as well as in \cite[p.~NP236]%
{MorGao:16}, where the indicated reference is \textquotedblleft Wu, C, Li, C,
and Gao, DY. Canonical primal-dual method for solving nonconvex minimization
problems. arXiv:1212.6492, 2012." Observe that the main difference between
arXiv:1212.6492 and reference \textquotedblleft\lbrack15]" of \cite{MorGao:17}
consists in the list of the authors, the content being practically the same.}

In fact, the same example given in \cite{VoiZal:11} but for $n=1$ shows that
even the results from \cite{MorGao:12}--\cite{MorGao:17} do not provide the
global solution of problem $(\mathcal{P})$. Indeed, as in \cite[p.~600]%
{VoiZal:11}, take $\gamma:=\sqrt{6}/96$; because $n=1$, we have that
$c=1\in\mathbb{R}$. Then the critical points of $\Xi$ with $(\lambda
,\mu,\varsigma)\in\mathcal{S}_{a}$ are, as indicated in \cite[p.~600]%
{VoiZal:11}, the following:%
\begin{gather*}
(\overline{x}_{1},\overline{y}_{1},\overline{\lambda}_{1},\overline{\mu}%
_{1},\overline{\varsigma}_{1}):=\big(1,1+\tfrac{1}{2}\sqrt{6},\tfrac{1}%
{2}\sqrt{6},\tfrac{48}{13},-\tfrac{1}{4}\big),\\
(\overline{x}_{2},\overline{y}_{2},\overline{\lambda}_{2},\overline{\mu}%
_{2},\overline{\varsigma}_{2}):=\big(-1,1+\tfrac{1}{2}\sqrt{6},-2-\tfrac{1}%
{2}\sqrt{6},\tfrac{16}{13}(3+2\sqrt{6}),-\tfrac{1}{4}\big)\\
(\overline{x}_{3},\overline{y}_{3},\overline{\lambda}_{3},\overline{\mu}%
_{3},\overline{\varsigma}_{3}):=\left(
1,2.603797322,1.603797322,-3.701\,325488,0.2860829239\right)  ,\\
(\overline{x}_{4},\overline{y}_{4},\overline{\lambda}_{4},\overline{\mu}%
_{4},\overline{\varsigma}_{4}):=\left(
-1,2.603797322,-3.603797322,-8.317027781,0.2860829239\right)  .
\end{gather*}

Using Corollary \ref{c-lkkt}, we have that $(\overline{\lambda}_{i}%
,\overline{\mu}_{i},\overline{\varsigma}_{i})$ with $i\in\overline{1,4}$ are
the only critical points of $D$ $(=\Pi^{d})$. For $i\in\{1,3\}$ we have that
$(1+\overline{\lambda}_{i})(1+\overline{\mu}_{i}\overline{\varsigma}_{i}%
)-1<0$, while for $i\in\{2,4\}$ we have that $1+\overline{\lambda}_{i}<0$ and
so $(\overline{\lambda},\overline{\mu},\overline{\varsigma})\notin
\mathcal{S}_{a}^{+}$ ($\mathcal{S}_{a}^{+}$ defined in \cite[Eq.~(18)]%
{MorGao:12} and \cite[Eq.~(19)]{MorGao:17}) and $(\overline{\lambda}%
,\overline{\mu},\overline{\varsigma})\notin\mathcal{S}_{c}^{+}$ ($\mathcal{S}%
_{c}^{+}$ defined in \cite[Eq.~(30)]{MorGao:16}). Therefore, the unique
solution $(1,1+\tfrac{1}{2}\sqrt{6})$ of problem $(\mathcal{P})$ is not
provided by either \cite[Th.~2.3]{MorGao:12}, or \cite[Th.~2]{MorGao:17}, or
\cite[Th.~3]{MorGao:16}. The use of the perturbation method suggested in these
papers is useless for this example.

\section{Conclusions}

-- We provided a rigorous treatment (study) for constrained minimization
problems using the Canonical duality theory developed by DY Gao.

-- Proposition \ref{p-perfdual} shows that the so-called perfect duality holds
under quite mild assumptions on the data of the problem; however, in our
opinion this formula is not very useful because for the found element
$(\overline{x},\overline{\lambda},\overline{\sigma})$, $\overline{x}$ could
not be feasible for the primal problem and/or $(\overline{\lambda}%
,\overline{\sigma})$ could not be feasible for the dual problem.

-- Proposition \ref{pei-i} and Remark \ref{rem-cdt} show that even if CDT can
be used for equality and/or inequality constrained optimization problems, it
is more appropriate for problems with inequality constraints.

-- The most important drawback of CDT is that it could find at most those
solutions of the primal problem for which all non quadratic constraints are
active; even more, the Lagrange multipliers corresponding to non quadratic
constraints must be strictly positive.

-- Moreover, the solutions found using CDT are among those found using the
usual Lagrange multipliers method. Using the \textquotedblleft extended
Lagrangian" $\Xi$ could be useful to decide if the found $\overline{x}$ is a
global minimizer of the primal problem.

-- The consideration of the dual function $D$ does not seem to be useful for
constrained minimization problems with at least one non quadratic constraint
because $D$ is not concave, unlike the case of quadratic constraints.

\medskip

\textbf{Acknowledgement} We thank Prof.\ Marius Durea for reading a previous
version of the paper and for his useful remarks.

\end{document}